\documentclass[a4paper,11pt]{article}
\usepackage{amsmath}
\usepackage{amsfonts}
\usepackage{amssymb}
\usepackage{latexsym}
\usepackage{epsfig}
\usepackage{graphicx}
\usepackage{oldgerm}
\usepackage{theorem}

\def\C{\mathbb{C}}
\def\R{\mathbb{R}}
\def\N{\mathbb{N}}
\def\Z{\mathbb{Z}}
\def\A{\mathbb{A}}

\def\cA{{\cal A}}

\def\x{\mib{x}}
\def\y{\mib{y}}
\def\u{\mib{u}}
\def\v{\mib{v}}

\def\f{\mib{f}}
\def\t{\mib{t}}
\def\V{\mib{V}}
\def\X{\mib{X}}

\def\bW{\mib{W}}
\def\bZ{\mib{Z}}

\def\s{\mib{s}}
\def\r{\mib{r}}

\def\rP{{\rm P}}
\def\rE{{\rm E}}
\def\cF{{\cal F}}

\def\P{{\mathbb P}}
\def\E{{\mathbb E}}

\def\1{{\bf 1}}

\def\mbK{\mathbb{K}}
\def\bE{{\bf E}}
\def\bP{{\bf P}}
\def\bK{{\bf K}}
\def\bG{{\bf G}}

\def\bbW{{\bf W}}

\def\mM{\mathfrak{M}}

\def\cM{{\cal M}}
\def\cD{{\cal D}}
\def\cG{{\cal G}}
\def\cK{{\cal K}}
\def\cS{{\cal S}}

\def\rC{{\rm C}}
\def\rS{{\rm S}}

\def\mX{\mathfrak{X}}

\def\chX{\check{X}}
\def\chW{\check{W}}
\def\n{\mib{n}}

\theorembodyfont{\itshape}

\newtheorem{thm}{Theorem}[section]
\newtheorem{lem}[thm]{Lemma}
\newtheorem{cor}[thm]{Corollary}
\newtheorem{prop}[thm]{Proposition}

\newtheorem{df}[thm]{Definition}

\newcommand{\mib}[1]{\mbox{\boldmath $#1$}}
\newcommand{\SSC}[1]{\section{#1}\setcounter{equation}{0}}
\newcommand{\qed}{\hbox{\rule[-2pt]{3pt}{6pt}}}

\begin{document}

\title{\bf 
Elliptic Determinantal Process of Type A
}
\author{
Makoto Katori
\footnote{
Department of Physics,
Faculty of Science and Engineering,
Chuo University, 
Kasuga, Bunkyo-ku, Tokyo 112-8551, Japan;
e-mail: katori@phys.chuo-u.ac.jp
}}
\date{29 September 2014}
\pagestyle{plain}
\maketitle
\begin{abstract}
We introduce an elliptic extension of Dyson's Brownian motion model,
which is a temporally inhomogeneous diffusion process of 
noncolliding particles defined on a circle.
Using elliptic determinant evaluations related to the reduced affine root system
of types $A$, we give determinantal martingale representation (DMR)
for the process, when it is started at the configuration 
with equidistant spacing on the circle. 
DMR proves that the process is determinantal
and the spatio-temporal correlation kernel is determined.
By taking temporally homogeneous limits of the present
elliptic determinantal process, trigonometric and hyperbolic
versions of noncolliding diffusion processes are studied.
\end{abstract}

Alcove of affine Weyl group


\SSC{Introduction \label{sec:introduction}}

Eigenvalue distributions of random-matrix ensembles provide
important examples of {\it determinantal point processes},
in which any correlation function is given by a determinant specified by
a single continuous function called the {\it correlation kernel}
\cite{Meh04,For10,Sos00,ST03,BKPV09}.
Dyson's Brownian motion model with parameter $\beta=2$ \cite{Dys62,Spo87}, 
which we simply call {\it the Dyson model} in this paper,  and 
other noncolliding diffusion processes \cite{Bru91,Gra99,KO01,Joh02,KT04} are
dynamical extensions of random-matrix ensembles. 
There any {\it spatio-temporal} correlation function is
expressed by determinant \cite{EM98,NF98} and such processes are said to
be {\it determinantal} \cite{BR05,KT10}.
The noncolliding diffusion processes have attracted much attention 
in probability  theory also by the fact that they are realized as
$h$-transforms in the sense of Doob of absorbing particle
systems in the Weyl chambers \cite{Gra99,KO01,KT04}. 
The relationship between the above mentioned integrability 
as spatio-temporal models and $h$-transform
constructions as stochastic processes 
has been clarified by introducing a notion of 
{\it determinantal martingales} in \cite{KT13,Kat13a,Kat13b}.
The purpose of the present paper is to report
{\it elliptic extensions} of these determinantal processes.
Since the Dyson model can be regarded as 
a multivariate extension of the three-dimensional Bessel process,
BES(3) \cite{KT11,BBO05}, 
first we discuss an elliptic extension of BES(3).

Let $i=\sqrt{-1}, v, \tau \in \C$ and put
\begin{equation}
z=z(v)=e^{\pi i v}, \quad q=q(\tau)=e^{\pi i \tau}.
\label{eqn:ztau1}
\end{equation}
The Jacobi theta function $\vartheta_1$ is defined as
\begin{eqnarray}
\vartheta_1(v; \tau)
&=& i \sum_{n \in \Z} (-1)^n q^{(n-(1/2))^2} z^{2n-1}
\nonumber\\
&=& 2 \sum_{n=1}^{\infty} (-1)^{n-1}
e^{\pi i \tau (n-(1/2))^2} \sin \{(2n-1) \pi v\}.
\label{eqn:Jacobi3}
\end{eqnarray}
(Note that the present function $\vartheta_1(v; \tau)$ is represented
as $\vartheta_1(\pi v,q)$ in \cite{WW27}.)
For $\Im \tau >0$, $\vartheta_1(v; \tau)$
is holomorphic for $|v| < \infty$
and satisfies the partial differential equation
\begin{equation}
\frac{\partial \vartheta_1(v; \tau)}{\partial \tau}=
\frac{1}{4 \pi i} \frac{\partial^2 \vartheta_1(v; \tau)}{\partial v^2}.
\label{eqn:Jacobi_eq}
\end{equation}
With parameters $N \in \N \equiv \{1,2, \dots\}, \alpha > 0$, and
$0 < t_{*} < \infty$, we introduce 
the following function of $(t, x) \in [0, t_{*}) \times \R$, 
\begin{eqnarray}
A_N^{\alpha}(t_{*}-t,x) &=& 
\left[ \frac{1}{2 \pi r} \frac{d}{dv}
\log \vartheta_1(v; \tau) \right]_{v=x/\alpha, \tau=2 \pi i N(t_{*}-t)/\alpha^2}
\nonumber\\
&=& \frac{1}{2 \pi r}
\frac{\vartheta_1'(x/\alpha ; 2 \pi i N(t_{*}-t)/\alpha^2)}
{\vartheta_1(x/\alpha ; 2 \pi i N(t_{*}-t)/\alpha^2)},
\label{eqn:A2}
\end{eqnarray}
where $\vartheta_1'(v; \tau)=d \vartheta_1(v; \tau)/dv$.
As a function of $x \in \R$, it is odd, 
\begin{equation}
A_N^{\alpha}(t_{*}-t,-x)=-A_N^{\alpha}(t_{*}-t, x),
\label{eqn:A_odd}
\end{equation}
and periodic with period $\alpha$
\begin{equation}
A_N^{\alpha}(t_{*}-t, x+m \alpha)=A_N^{\alpha}(t_{*}-t, x),
\quad m \in \Z.
\label{eqn:A_periodic}
\end{equation}
It has only simple poles at
$x= m \alpha, m \in \Z$, 
and simple zeroes at
$x=(m+1/2) \alpha, m \in \Z$.

Let $r > 0$. Suppose that $\chX(t), t \in [0, t_{*})$ satisfies
the stochastic differential equation (SDE) 
\begin{equation}
d \chX(t)=dB(t)+ A_1^{2 \pi r}(t_{*}-t, \chX(t)) dt
\label{eqn:eBES1}
\end{equation}
started at $\chX(0)=x \in (0, 2 \pi r)$,
where $B$ denotes the one-dimensional standard
Brownian motion (BM). 
(From now on, BM means a one-dimensional standard
Brownian motion unless specially mentioned.) 
By periodicity (\ref{eqn:A_periodic}) with period $\alpha=2 \pi r$, 
$\chX$ can be considered to describe 
a diffusion process of a particle moving around a circle with radius $r>0$;
$\rS^1(r)=\{x \in \R : x + 2 \pi r = x\}$. 
Note that this system is temporally inhomogeneous
defined only in a time interval $[0, t_{*})$.
Independently of $t \in [0, t_{*})$ and $N$, however,
we have 
\begin{eqnarray}
&& A_N^{2 \pi r}(t_{*}-t,x) \sim \frac{1}{x} \quad \mbox{as $x \downarrow 0$},
\nonumber\\
&& A_N^{2 \pi r}(t_{*}-t,x) \sim - \frac{1}{2 \pi r -x}
\quad \mbox{as $x \uparrow 2 \pi r$}.
\label{eqn:A_asym}
\end{eqnarray} 
It implies that the behavior of $\chX \in (0, 2 \pi r)$ in the vicinity of 0
(and $2 \pi r$) 
is similar to that of BES(3) near 0.
We define a process $X \in [0, 2 \pi r)$ by
\begin{equation}
X(t)=\chX(t) \quad \mbox{mod $2 \pi r$}, \quad t \in [0, t_*).
\label{eqn:Def_X1}
\end{equation}
It gives a Markov process showing a position on 
the circumference $[0, 2 \pi r)$ of $\rS^1(r)$.
We write the probability law of $X(t), t \in [0, t_{*})$ 
started at $x=X(0) \in (0, 2 \pi r)$ as $\P_x$.

The backward Kolmogorov equation
for the SDE (\ref{eqn:eBES1}) is given as
\begin{equation}
- \frac{\partial u(s, x)}{\partial s}
= \frac{1}{2} \frac{\partial^2 u(s, x)}{\partial x^2}
+ A_1^{2 \pi r}(t_{*}-s, x) \frac{\partial u(s, x)}{\partial x},
\quad 0 \leq s < t_{*},
\quad x \in (0, 2 \pi r).
\label{eqn:Kol1}
\end{equation}
Let $q(t-s, y|x)$ be a solution of diffusion equation
$-\partial v(s,x)/\partial s=(1/2) \partial^2 v(s,x)/\partial x^2$,
$0 \leq s \leq t < \infty, x \in (0, 2 \pi r)$ satisfying
$\lim_{s \uparrow t} v(s,x)=\delta_y(\{x\})$, $y \in (0, 2 \pi r)$.
Then
\begin{equation}
u(s, x) = p(t, y|s,x) \equiv q(t-s, y|x)
\frac{\vartheta_1(y/2 \pi r; i (t_{*}-t)/2 \pi r^2)}
{\vartheta_1(x/2 \pi r; i (t_{*}-s)/2 \pi r^2)}, \quad 0 \leq s \leq t < t_{*}, 
\label{eqn:Kol2}
\end{equation}
solves (\ref{eqn:Kol1}) under the condition
$\lim_{s \uparrow t} u(s,x)=\delta_y(\{x\})$, $y \in (0, 2 \pi r)$,
since $\vartheta_1(v; \tau)$ satisfies (\ref{eqn:Jacobi_eq}).
For $0 < y < 2 \pi r$,
$\vartheta_1(y/2 \pi r; i(t_{*}-t)/2 \pi r^2) >0$
and it has simple zeroes at $y=0$ and $y= 2 \pi r$,
$t \in [0, t_{*})$.
Then, if $q(t-s, y|x)$ is chosen as 
the transition probability density (tpd) of the absorbing BM
in the interval $[0, 2 \pi r]$ with absorbing walls
at $x=0$ and $x=2 \pi r$, 
(\ref{eqn:Kol2}) is strictly positive and finite for any $x, y \in (0, 2 \pi r)$,
$0 \leq s \leq t < t_{*}$,
and thus $p(t, y|s, x)$ gives the tpd for
the process $X(t), t \in [0, t_{*})$.
Let $W(t), t \geq 0$ be BM started at $x=W(0) \in (0, 2 \pi r)$,
where its probability law is denoted by $\rP_x$.
Consider a filtration $\{\cF_W(t): t \geq 0\}$ 
generated by $W(t), t \geq 0$, 
which satisfies the usual conditions, 
and introduce a stopping time
$T_W = \inf \{ t >0: W(t) \in \{0, 2 \pi r\} \}$.
Then the above fact implies that, for $t \in [0, t_{*})$, 
\begin{equation}
\P_x(X(t) \in dy)
= \rP_x(T_W > t, W(t) \in dy) 
\frac{\vartheta_1(y/2 \pi r; i (t_{*}-t)/2 \pi r^2)}
{\vartheta_1(x/2 \pi r; it_{*}/2 \pi r^2)},
\quad x, y \in (0, 2 \pi r).
\label{eqn:Doob1}
\end{equation}

Note that
$$
\vartheta_1(v; \tau) \sim 2 e^{\pi i \tau/4} \sin (\pi v)
\quad \mbox{as $\Im \tau \to +\infty$ ({\it i.e.}, $q=q(\tau)=e^{\pi i \tau} \to 0$)}.
$$
Thus, in the limit $t_{*} \to \infty$, 
(\ref{eqn:eBES1}) becomes 
a temporally homogeneous SDE, 
\begin{equation}
d \chX(t)=dB(t)+\frac{1}{2r} \cot \left( 
\frac{\chX(t)}{2r} \right) dt,
\quad t \geq 0, 
\label{eqn:cot1}
\end{equation}
and (\ref{eqn:Doob1}) becomes
\begin{equation}
\P_x(X(t) \in dy)
= \rP_x(T_W > t, W(t) \in dy) 
\frac{\sin(y/2 r)}
{\sin(x/2 r)},
\quad  t \in [0, \infty), \quad x, y \in (0, 2 \pi r).
\label{eqn:Doob2}
\end{equation}
If we take the further limit
$r \to \infty$ in (\ref{eqn:cot1}), we have
$X(t)=\chX(t), t \geq 0$ and 
\begin{equation}
dX(t)=dB(t)+\frac{dt}{X(t)},
\quad t \geq 0,
\label{eqn:BES1}
\end{equation}
which is the SDE for BES(3)
on $\R_+=\{x > 0 : x \in \R\}$, and
(\ref{eqn:Doob2}) becomes
\begin{equation}
\P_x(X(t) \in dy)
= \rP_x(T'_W > t, W(t) \in dy) 
\frac{y}{x},
\quad t \in [0, \infty), 
\quad x, y \in \R_+,
\label{eqn:Doob3}
\end{equation}
where $T'_W=\inf \{t >0 : W(t)=0\}$.
The relation (\ref{eqn:Doob3}) states that BES(3)
is the {\it Doob $h$-transform of the absorbing BM}
in $[0,  \infty)$ with an absorbing wall at the origin,
where the harmonic function is given by
$h(x)=x, x \geq 0$.
We regard (\ref{eqn:cot1}) as a trigonometric extension,
and (\ref{eqn:eBES1}) as an elliptic extension
of (\ref{eqn:BES1}), respectively.
The equality (\ref{eqn:Doob3}) is generalized to
(\ref{eqn:Doob2}) and (\ref{eqn:Doob1}), respectively.
We can also discuss a scaling limit realizing $q \to 1$,
in which hyperbolic version of (\ref{eqn:eBES1})
is obtained (see Section \ref{sec:sinq}).

Let $N \in \{2,3, \dots\}$ and consider the following bounded region
$$
\cA^{A_{N-1}}_{2 \pi r}=\{ \x=(x_1, \dots, x_N) \in \R^N :
x_1< x_2 < \cdots < x_N < x_1+2 \pi r \},
$$
which is called a {\it scaled alcove of the affine Weyl group of type}
$A_{N-1}$ (with scale $2 \pi r$) \cite{Gra02,Kra07}.
With an additional condition $x_1 \geq 0$, we also consider the space
\begin{eqnarray}
\cA_{[0, 2 \pi r)^N}
&=& \cA^{A_{N-1}}_{2 \pi r} \cap \{\x \in \R^N : x_1 \geq 0 \}
\nonumber\\
&=& \{ \x \in \R^N : 0 \leq x_1 < x_2 < \cdots < x_N < 2 \pi r \}.
\nonumber
\end{eqnarray}
Note that it is different from the scaled alcove of type $C_N$ 
defined by
$\cA^{C_N}_{2 \pi r}=\{\x \in \R^N: 0 < x_1 < x_2 < \cdots < x_N < 2 \pi r\}$ 
\cite{Gra02,Kra07}
which excludes $x_1=0$ from $\cA_{[0, 2 \pi r)^N}$.

Now we introduce an $N$-particle extension
of the above process, 
$\check{\X}^A(t)=(\chX^A_1(t), \dots, \chX^A_N(t))$, $t \in [0, t_{*})$.
Assume that the initial configuration is chosen in the alcove
$$
\check{\X}^A(0) =\u \in \cA^{A_{N-1}}_{2 \pi r},
$$
and an index $\delta \in \pi r \Z$ is determined so that
$$
\overline{u}_{\delta} \equiv \delta+\sum_{j=1}^N u_j \in (0, 2 \pi r).
$$
Let 
$$
\overline{X}^A_{\delta}(t)=\delta+ \sum_{j=1}^N \chX^A_j(t), \quad
t \in [0, t_{*}).
$$
Then $\check{\X}^A(t), t \in [0, t_{*})$ is defined
as a solution of the following set of SDEs on $\R$,
\begin{equation}
d \chX^A_j(t) = dB_j(t)+ \sum_{\substack{1 \leq k \leq N, \cr k \not=j}}
A_N^{2 \pi r} (t_{*}-t, \chX^A_j(t)- \chX^A_k(t)) dt
+ A_N^{2 \pi r}(t_{*}-t, \overline{X}^A_{\delta}(t)) dt, 
\label{eqn:SDEe1}
\end{equation}
$1 \leq j \leq N, t \in [0, t_{*})$, 
where $B_j, 1 \leq j \leq N$ 
are independent BMs on $\R$.
By (\ref{eqn:A_odd}), 
$\sum_{1 \leq j, k \leq N, j \not= k}A_N^{2 \pi r}(t_{*}-s, x_j-x_k)=0$, 
and the summation of  (\ref{eqn:SDEe1}) over $j=1,2, \dots, N$ gives
\begin{equation}
d \overline{X}^A_{\delta}(t)
= \sqrt{N} dB(t) + N A_N^{2 \pi r}(t_{*}-t, \overline{X}^A_{\delta}(t)) d t,
\quad t \in [0, t_{*}), 
\label{eqn:SDEe2}
\end{equation}
where $B$ is BM on $\R$.
We then define the process
$\X^A(t)=(X^A_1(t), \dots, X^A_N(t)) \in [0, 2 \pi r)^N, t \in [0, t_*)$ by
\begin{equation}
X^A_j(t)=\chX^A_j(t) \quad \mbox{mod $2 \pi r$}, \quad
1\leq j \leq N, \quad t \in [0, t_*).
\label{eqn:Def_Xjs}
\end{equation}
It represents a Markov process showing the positions of $N$ particles
on the circumference $[0, 2 \pi r)$ of $\rS^1(r)$.

Let $\mM([0, 2 \pi r))$ be the space of nonnegative integer-valued Radon measures 
on the interval $[0, 2 \pi r)$, 
which is a Polish space with the vague topology.
Any element $\xi$ of $\mM([0, 2 \pi r))$ can be represented as
$\xi(\cdot) = \sum_{j \geq 1}\delta_{x_j}(\cdot)$, in which
the sequence of points in $[0, 2 \pi r)$, $\x =(x_j)_{j \geq 1}$, 
satisfies $\xi(K)=\sharp\{x_j : x_j \in K \} < \infty$ 
for any subset $K \subset [0, 2 \pi r)$.
Now we consider the process $\X^A(t)$
as $\mM([0, 2 \pi r))$-valued process
and write it as
\begin{equation}
\Xi^A(t, \cdot)=\sum_{j=1}^N \delta_{X^A_j(t)}(\cdot),
\quad t \in [0, t_{*}).
\label{eqn:Xi1}
\end{equation}
The probability law of $\Xi^A(t, \cdot), t \in [0, t_{*})$
starting from a fixed configuration $\xi \in \mM([0, 2 \pi r))$
is denoted by $\P^A_{\xi}$ 
and the process
specified by the initial configuration
is expressed by
$(\Xi^A(t), t \in [0,t_{*}), \P^A_{\xi})$.
The expectations with respect to $\P^A_{\xi}$
is denoted by $\E^A_{\xi}$.
We introduce a filtration $\{\cF_{\Xi^A}(t) : t \in [0,t_{*}) \}$
generated by $\Xi^A(t), t \in [0, t_*)$, 
which satisfies the usual conditions. 
Let $\rC([0, 2 \pi r))$ be the set of all continuous
real-valued functions on $[0, 2 \pi r)$.
We set 
$$
\mM_0([0, 2 \pi r))= \{ \xi \in \mM([0, 2 \pi r)) : 
\xi(\{x\}) \leq 1 \mbox { for any }  x \in [0, 2 \pi r) \},
$$
which denotes a collection of configurations
without any multiple points.

In the present paper, we study the case that
the initial state
$\eta=\sum_{j=1}^N \delta_{v_j}$ 
is corresponding to 
the configuration with equidistant spacing on $\rS^1(r)$;
\begin{equation}
v_j=\frac{2 \pi r}{N}(j-1), \quad
1 \leq j \leq N, 
\label{eqn:ed1}
\end{equation}
and we will prove that the process
$(\Xi^A(t), t \in [0, t_{*}), \P^A_{\eta})$ is determinantal. 
In this case, $\sum_{j=1}^N v_j=\pi r(N-1)$
and the index $\delta \in \pi r \Z$ is determined as
\begin{equation}
\delta=- \pi r (N-2), 
\label{eqn:delta_A}
\end{equation}
so that $\overline{v}_{\delta}= \pi r \in (0, 2 \pi r)$.
We will present that this determinantal process
can be considered as an elliptic extension of the Dyson model.
The key lemmas to construct the elliptic determinantal process
are obtained from the {\it elliptic determinant evaluations}
related to infinite families of irreducible reduced affine root systems
studied in \cite{For90a,For90b,TV97,War02,RS06,Kra05,For10}.
According to Macdonald's classification of
reduced affine root systems \cite{Mac72},
the present process is related to the system of type $A_{N-1}$.
Further study concerning other types is in progress.
Connection between the elliptic determinantal processes
and probabilistic discrete models
with elliptic weights \cite{Sch07,BGR10,Betea11} will be an
interesting future problem.
We note that the function $A_N^{\alpha}(t_{*}-t, z)$ can be regarded
as Villat's kernel for an annulus $\A_q=\{z \in \C : q < |z| < 1\}$
with $0 < q=e^{-2 \pi^2 N(t_*-t)/\alpha^2} < 1$.
It is the reason why it also appears in the study of
stochastic Komatu-Loewner evolution in doubly connected domains
\cite{Zhan04,BF08}.

The paper is organized as follows.
In Section \ref{sec:preliminaries} 
preliminaries of elliptic functions
and their related functions are given. Useful determinantal
identities are obtained from the elliptic determinant evaluations
related to the affine root system of types $A_{N-1}$ 
\cite{For90a,For90b,TV97,War02,RS06,Kra05,For10}.
There the generalized $h$-transform and determinantal martingale
for $(\Xi^A(t), t \in [0, t_{*}), \P^A_{\eta})$ are introduced.
The main results are given in Section \ref{sec:results}.
First the determinantal martingale representation (DMR) is given
for $(\Xi^A(t), t \in [0, t_{*}), \P^A_{\eta})$
(Theorem \ref{thm:DMR}). 
As a result of Theorem 1.3 of \cite{Kat13a}, 
DMR proves that
the process is determinantal and the correlation kernel is determined
(Corollary \ref{thm:kernel}).
Explicit expression for the correlation kernel
of $(\Xi^A(t), t \in [0, t_*), \P^A_{\eta})$ is shown
and the infinite-particle limit is discussed.
In Section \ref{sec:reductions} 
the temporally homogeneous limit $t_{*} \to \infty$
is studied both in the level of SDE
and in the level of determinantal process.
We study the system of noncolliding Brownian motions on a circle,
$\widehat{\Xi}^A(t), t \in [0, \infty)$, 
obtained from $\Xi^A(t), t \in [0, t_*)$ by this reduction.
It is different from the dynamical CUE (circular unitary ensemble) model
studied in \cite{HW96,NF03,Kat13a}.
Finally in Section \ref{sec:sinq}, the results are expressed
by using Gosper's $q$-sine function \cite{Gos01}
as well as the $q$-gamma function \cite{AAR99}.
Then $q \to 1$ limit is discussed, in which temporally homogeneous
processes expressed by hyperbolic functions are obtained.

\SSC{Preliminaries \label{sec:preliminaries}}
\subsection{Elliptic functions and their related functions \label{sec:elliptic}}
The Jacobi theta function $\vartheta_1$ defined by
(\ref{eqn:Jacobi3}) has the following infinite-product expressions, 
\begin{eqnarray}
\vartheta_1(v ; \tau)
&=& -i q^{1/4} q_0 z \prod_{j=1}^{\infty} 
(1-q^{2j} z^2) (1-q^{2j-2}/z^2)
\nonumber\\
&=& 2 q^{1/4} q_0 \sin(\pi v)
\prod_{j=1}^{\infty}
(1-2 q^{2j} \cos(2 \pi v) + q^{4j})
\label{eqn:theta_p1}
\end{eqnarray}
with 
$$
q_0= q_0(\tau) \equiv \prod_{n=1}^{\infty} (1-q^{2n}), \quad q=e^{\pi i \tau}.
$$
It is easy to see from these expressions that,
when $\Re \tau=0, \Im \tau >0$ (that is, $0 < q < 1$), 
$\vartheta_1(v; \tau)>0$ for $v \in (0, 1)$
and it has simple zeroes at $v=0$ and $v=1$.
It is odd with respect to $v$
\begin{equation}
\vartheta_1(-v; \tau)=-\vartheta_1(v; \tau), 
\label{eqn:Jacobi4}
\end{equation}
and has quasi-periodicity
\begin{eqnarray}
\label{eqn:Jacobi5a}
\vartheta_1(v+1; \tau)
&=& -\vartheta_1(v; \tau),
\\
\vartheta_1(v+\tau; \tau)
&=&-\frac{1}{z^2 q} \vartheta_1(v; \tau)
=-e^{-\pi i (2v+\tau)} \vartheta_1(v; \tau).
\nonumber
\end{eqnarray}
We define 
\begin{eqnarray}
\vartheta_0(v; \tau) &\equiv& -i e^{\pi i(v+\tau/4)} \vartheta_1 \left(v + \frac{\tau}{2}; \tau \right)
= \sum_{n \in \Z} (-1)^n q^{n^2} z^{2n},
\nonumber\\
\vartheta_2(v; \tau) &\equiv& \vartheta_1 \left(v+ \frac{1}{2}; \tau \right)
= \sum_{n \in \Z} q^{(n-(1/2))^2} z^{2n-1},
\nonumber\\
\vartheta_3(v; \tau) &\equiv& e^{\pi i(v+\tau/4)} \vartheta_1 \left(v + \frac{1+\tau}{2}; \tau \right)
= \sum_{n \in \Z} q^{n^2} z^{2n},
\label{eqn:theta_0_2_3}
\end{eqnarray}
as usual.
(Note that the present functions $\vartheta_{\mu}(v; \tau)$, $\mu=1,2,3$ are denoted 
by $\vartheta_{\mu}(\pi v,q)$, and $\vartheta_0(v; \tau)$ by $\vartheta_4(\pi v, q)$ 
in \cite{WW27}.)
They solve the partial differential equation
\begin{equation}
\frac{\partial \vartheta_{\mu}(v; \tau)}{\partial \tau}=
\frac{1}{4 \pi i} \frac{\partial^2 \vartheta_{\mu}(v; \tau)}{\partial v^2},
\quad \mu=0,1,2,3.
\label{eqn:Jacobi_equation}
\end{equation}

We will use the following formulas; $n \in \N$,
$q=e^{\pi i \tau}$, 
\begin{eqnarray}
\label{eqn:formula1}
&& \prod_{j=0}^{n-1} \vartheta_1 (v+j/n; \tau)
=\frac{q_0^n}{(q^n)_0} \vartheta_1(nv; n \tau),
\\
\label{eqn:formula2}
&& 
\prod_{j=1}^{n-1} \vartheta_1(j/n; \tau)
= \frac{n q_0^n}{(q^n)_0}
\frac{\vartheta_1'(0; n \tau)}{\vartheta_1'(0; \tau)},
\end{eqnarray}
where
\begin{equation}
\vartheta_1'(0; \tau)
\equiv \left. \frac{\partial \vartheta_1(v; \tau)}{\partial v}
\right|_{v=0}
=2 \pi \sum_{j=1}^{\infty} (-1)^{j-1} (2j-1) q^{n-(1/2))^2},
\label{eqn:theta_dash}
\end{equation}
and
\begin{equation}
\frac{\vartheta_1(v+w; \tau) \vartheta_1'(0; \tau)}
{\pi \vartheta_1(v; \tau) \vartheta_1(w; \tau)}
= \cot(\pi v) + \cot(\pi w)
+ 4 \sum_{\ell=1}^{\infty} \sum_{m=1}^{\infty} q^{2 \ell m}
\sin [ 2 \pi (\ell v+ m w)].
\label{eqn:formula3}
\end{equation}
The formula (\ref{eqn:formula1}) is obtained from
Eq.(2.3) of \cite{RS06} (see (\ref{eqn:root1}) below)
and (\ref{eqn:formula2}) is obtained as its $v \to 0$
limit, since $\vartheta_1(v; \tau)$ has a simple root at $v=0$.
The formula (\ref{eqn:formula3}) is found in
`Miscellaneous Examples' of Chapter XXI in \cite{WW27}.

Let $\omega_1$ and $\omega_3$ be fundamental periods and set
\begin{eqnarray}
\tau &=& \frac{\omega_3}{\omega_1},\quad \Im \tau >0,
\nonumber\\
\Omega_{m,n} &=& 2 m \omega_1+2n \omega_3,
\quad m,n \in \Z.
\nonumber
\end{eqnarray}
The Weierstrass $\wp$ function and zeta function $\zeta$ 
are defined as the following meromorphic functions
\begin{eqnarray}
\wp(z) &=& \wp(z| 2 \omega_1, 2 \omega_3)
\nonumber\\
&=& \frac{1}{z^2}
+\sum_{(m,n) \in \Z^2 \setminus \{(0,0)\}}
\left[ \frac{1}{(z-\Omega_{m,n})^2}-\frac{1}{{\Omega_{m,n}}^2} \right],
\nonumber\\
\zeta(z) &=& \zeta(z| 2 \omega_1, 2 \omega_3)
\nonumber\\
&=& \frac{1}{z}
+\sum_{(m,n) \in \Z^2 \setminus \{(0,0)\}}
\left[ \frac{1}{z-\Omega_{m,n}}+\frac{1}{\Omega_{m,n}}
+\frac{z}{{\Omega_{m,n}}^2} \right].
\nonumber
\end{eqnarray}
Let
$\omega_2=-(\omega_1+\omega_3)$,
and put
$\zeta(\omega_{\nu})=\eta_{\nu}, \nu=1,2,3$.
The relation
$\eta_1+\eta_2+\eta_3=0$
holds.
By definition, $\wp(z)$ is an elliptic function
with fundamental periods $\omega_1, \omega_2$ and $\omega_3$,
$$
\wp(z+2 \omega_{\nu})=\wp(z), \quad
\nu=1,2,3,
$$
and it is even,
$\wp(-z)=\wp(z)$.
The function $\zeta$ is odd,
$\zeta(-z)=-\zeta(z)$,  and
is quasi-periodic in the sense
$$
\zeta(z+2 \omega_{\nu})=\zeta(z)+2 \eta_{\nu}, \quad
\nu=1,2,3.
$$
By definition, 
\begin{equation}
\wp(z)=-\zeta'(z). 
\label{eqn:wp2}
\end{equation}
Moreover, the relation
\begin{equation}
\{\zeta(z+u)-\zeta(z)-\zeta(u)\}^2
=\wp(z+u)+\wp(z)+\wp(u)
\label{eqn:zeta_P0}
\end{equation}
holds (see Section 20.41 in \cite{WW27}). 
From this, we obtain the following identity.
\begin{lem}
\label{thm:zeta_P1}
For $a, b, c \in \C$, 
\begin{eqnarray}
&& \zeta(a-b)\zeta(a-c)+\zeta(b-a)\zeta(b-c)+\zeta(c-a)\zeta(c-b)
\nonumber\\
&& \quad = \frac{1}{2} \Big\{
\zeta(a-b)^2+\zeta(b-c)^2+\zeta(a-c)^2 \Big\}
-\frac{1}{2} \Big\{ \wp(a-b)+\wp(b-c)+\wp(a-c) \Big\}.
\nonumber
\end{eqnarray}
\end{lem}
\noindent{\it Proof.} 
Put $z=a-b, u=b-c$ in (\ref{eqn:zeta_P0}).
Then by the fact that $\zeta$ is odd, the equality is
derived. \qed
\vskip 0.3cm
The following relation is established,
$$
\zeta(z| 2\omega_1, 2 \omega_3) - \frac{\eta_1 z}{\omega_1}
= \frac{1}{2 \omega_1} 
\left. \frac{d}{dv} \log \vartheta_1(v; \tau) \right|_{v=z/2\omega_1},
\quad \tau=\frac{\omega_3}{\omega_1}.
$$
It has the following expansion with respect to $q = q(\tau) = e^{\pi i \tau}$, 
\begin{equation}
\zeta(z| 2\omega_1, 2 \omega_3) - \frac{\eta_1 z}{\omega_1}
= \frac{\pi}{2 \omega_1} \cot \left( \frac{\pi z}{2 \omega_1} \right)
+ \frac{2 \pi}{\omega_1}
\sum_{n=1}^{\infty} \frac{q^{2n}}{1-q^{2n}} \sin \left( \frac{n \pi z}{\omega_1} \right).
\label{eqn:relation2}
\end{equation}
The function $A_N^{\alpha}(t_{*}-t, x)$ defined by (\ref{eqn:A2}) is then expressed as
\begin{eqnarray}
A_N^{\alpha}(t_{*}-t,x) &=&  \left[ \zeta(x|2 \omega_1, 2 \omega_3) 
-\frac{\eta_1 x}{\omega_1} \right]_{\omega_1=\alpha/2, \omega_3=\pi i N(t_{*}-t)/\alpha}
\nonumber\\
&=& \zeta\left(x \left| \alpha, \frac{2 \pi i N(t_{*}-t)}{\alpha} \right. \right) 
- \frac{2 \eta_1(t_{*}-t) x}{\alpha},
\label{eqn:A1}
\end{eqnarray}
for $t \in [0, t_{*})$, $x \in \R$, 
where
\begin{eqnarray}
\eta_1(t_{*}-t)=\eta_1(t_{*}-t; N, \alpha)
&=& \frac{\pi^2}{\omega_1} \left. \left( \frac{1}{12}-2 \sum_{n=1}^{\infty} \frac{n q^{2n}}{1-q^{2n}} \right)
\right|_{\omega_1=\alpha/2, q=e^{-2 \pi^2 N(t_{*}-t)/\alpha^2}}
\nonumber\\
&=& \frac{2 \pi^2}{\alpha} \left( \frac{1}{12}
- 2 \sum_{n=1}^{\infty} \frac{n e^{-4 \pi^2 n N(t_{*}-t)/\alpha^2}}{1-e^{-4\pi^2 n N(t_{*}-t)/\alpha^2}}\right).
\label{eqn:eta1b}
\end{eqnarray}
The formula (\ref{eqn:relation2}) gives
\begin{equation}
A_N^{\alpha}(t_{*}-t,x) 
= \frac{\pi}{\alpha} \cot \left( \frac{\pi x}{\alpha} \right)
+ \frac{4 \pi}{\alpha} \sum_{n=1}^{\infty} \frac{e^{- 4 \pi^2 n N(t_{*}-t)/\alpha^2}}{1-e^{-4 \pi^2 n N (t_{*}-t)/\alpha^2}}
\sin \left( \frac{2 \pi n x}{\alpha} \right),
\quad t \in [0, t_{*}), \quad x \in \R.
\label{eqn:A3}
\end{equation}
From this expression, we can readily observe (\ref{eqn:A_odd}) and 
(\ref{eqn:A_periodic}) and other properties of $A_N^{\alpha}$.
In particular, we can see (\ref{eqn:A_asym}).

\subsection{Elliptic determinant identities \label{sec:det_identity}}

Let $p$ be a fixed complex 
number such that $0<|p|<1$.
We use the standard notations
\begin{eqnarray}
(a;p)_{\infty} &=& \prod_{j=0}^{\infty} (1-a p^j),
\nonumber\\
(a_1,\dots, a_n;p)_{\infty}
&=& (a_1;p)_{\infty} \cdots (a_n;p)_{\infty}.
\nonumber
\end{eqnarray}
Following \cite{War02,RS06}, here we use `multiplicative notation'
for theta functions,
\begin{eqnarray}
E(s; p) &=& (s, p/s; p)_{\infty},
\nonumber\\
E(s_1, \dots, s_n; p) &=& \prod_{j=1}^n E(s_j;p),
\nonumber\\
E(s_1 s_2^{\pm}; p) &=& E(s_1 s_2; p) E(s_1/s_2; p).
\nonumber
\end{eqnarray}
The function $E(s;p)$ is holomorphic for $s \not=0$.
The zero set of $E(s;p)$ is given by
$p^{\Z} \equiv \{p^j: j \in \Z\}$, 
and all zeroes are single.
The inversion formula
\begin{equation}
E(1/s; p)=-\frac{1}{s} E(s; p),
\label{eqn:inversion_E}
\end{equation} 
the quasi-periodicity
$E(ps; p)=-E(s; p)/s$,
and the Laurent expansion
$$
E(s; p)=\frac{1}{(p;p)_{\infty}}
\sum_{n \in \Z} (-1)^n p^{{n \choose 2}} s^n
$$
are known.
If $g_n$ denotes a primitive $n$th root of unity, $n \in \N$, 
the following equality holds \cite{RS06},
\begin{equation}
E(s^n; p^n)=\prod_{j=0}^{n-1} E(s g_n^j; p).
\label{eqn:root1}
\end{equation}

With (\ref{eqn:ztau1}), the function $E(s;p)$ is
related to the Jacobi theta function $\vartheta_1$ by
\begin{equation}
\vartheta_1(v; \tau)
= i q^{1/4} q_0 \frac{1}{z} E(z^2;q^2).
\label{eqn:JacobiB1}
\end{equation}
The equality (\ref{eqn:root1}) is then rewritten as 
(\ref{eqn:formula1}).

For $N \in \N$, 
$\s=(s_1, \dots, s_N) \in \C^N$, put
$$
W_{A_{N-1}}(\s; p)
= \prod_{1 \leq j < k \leq N} 
s_k E(s_j/s_k; p)
= (-1)^{N(N-1)/2} \prod_{1 \leq j < k \leq N}
s_j E(s_k/s_j;p),
$$
where the second equality is proved by
the inversion formula (\ref{eqn:inversion_E}). 
It is the {\it Macdonald denominator} for 
the reduced affine root systems of type $A_{N-1}$ \cite{RS06}.

We start with the following two lemmas, which are readily
obtained from the results given in
\cite{TV97,War02,RS06,Kra05}.

\begin{lem}
\label{thm:equality1}
For $\r=(r_1, \dots, r_N) \in \C^N$ and $\kappa \in \C$ assume that
\begin{equation}
\frac{r_j}{r_k} \notin p^{\Z}, \quad 1 \leq j \not= k \leq N,
\quad \mbox{and} \quad 
\kappa \prod_{j=1}^N r_j \notin p^{\Z}.
\label{eqn:ass_u1}
\end{equation}
Then for $\s \in \C^N$
\begin{eqnarray}
&& 
\det_{1 \leq j, k \leq N} 
\left[
\frac{E \left( \kappa s_j \prod_{m=1, m \not= k}^N r_m ;p \right)}
{E \left( \kappa \prod_{m=1}^N r_m ;p \right)}
\prod_{1 \leq \ell \leq N, \ell \not=k}
\frac{E(s_j/r_{\ell};p)}{E(r_k/r_{\ell};p)} 
\right]
\nonumber\\
&& \qquad
= \frac{E \left( \kappa \prod_{j=1}^N s_j ;p \right)
W_{A_{N-1}}(\s; p)}
{E \left( \kappa \prod_{j=1}^N r_j ;p \right)
W_{A_{N-1}}(\r; p)}.
\label{eqn:equality1}
\end{eqnarray}
\end{lem}
\noindent{\it Proof of Lemma \ref{thm:equality1}.} 
Let $s_1, \dots, s_N$, $a_1, \dots, a_{N-1}$, 
$b_1, \dots, b_N$, $c_2, \dots, c_N$,
and $\kappa$ be indeterminates, which satisfy
\begin{equation}
\prod_{k=1}^{j-1} a_k \cdot b_j \cdot
\prod_{\ell=j+1}^N c_{\ell}=\kappa, \quad
1 \leq j \leq N.
\label{eqn:TV1}
\end{equation}
The following equality holds \cite{TV97}
(see Corollary 4.5 and Remark 4.6 in \cite{RS06},
and see also Section 5.11 in \cite{Kra05}), 
\begin{eqnarray}
&& \det_{1 \leq j, k \leq N}
\left[ \prod_{\ell=1}^{k-1} E(a_{\ell} s_j; p) \cdot
E(b_k s_j; p) \cdot
\prod_{m=k+1}^N E(c_m s_j; p) \right]
\nonumber\\
&& \qquad
= E \left( \kappa \prod_{j=1}^N s_j ;p \right)
\prod_{k=2}^N E(b_k/c_k; p)
\prod_{1 \leq \ell < m \leq N} c_m s_m
E(s_{\ell}/s_m, a_{\ell}/c_m; p).
\label{eqn:TV2}
\end{eqnarray}
We put $a_j=1/r_j, 1 \leq j \leq N-1$,
$b_j=\alpha/r_j, 1 \leq j \leq N$,
$c_j=1/r_j, 2 \leq j \leq N$ with $\alpha \in \C$.
Then the condition (\ref{eqn:TV1}) gives
\begin{equation}
\prod_{j=1}^N r_j=\frac{\alpha}{\kappa}
\label{eqn:TVa1}
\end{equation}
and (\ref{eqn:TV2}) becomes
\begin{eqnarray}
&& \det_{1 \leq j, k \leq N} \left[
E(\alpha s_j/r_k; p) \prod_{1 \leq \ell \leq N, \ell \not=k}
E(s_j/r_{\ell}; p) \right]
\nonumber\\
&& \qquad 
=E \left( \alpha
\prod_{j=1}^N s_j/r_j ; p\right)
E(\alpha; p)^{N-1}
\prod_{1 \leq \ell < m \leq N}
\frac{s_m}{r_m}
E \left( \frac{s_{\ell}}{s_m}; p \right)
E \left( \frac{r_m}{r_{\ell}}; p \right).
\label{eqn:TVa2}
\end{eqnarray}
Under the assumption (\ref{eqn:ass_u1})
we divide the both sides of (\ref{eqn:TVa2}) by
$\displaystyle{E(\alpha; p)^N \prod_{1 \leq k, \ell \leq N, k \not= \ell}
E(r_k/r_{\ell}; p)}$ 
and obtain the equality
$$
\det_{1 \leq j, k \leq N}
\left[ \frac{E(\alpha s_j/r_k; p)}{E(\alpha; p)}
\prod_{1 \leq \ell \leq N, \ell \not=k}
\frac{E(s_j/r_{\ell}; p)}{E(r_k/r_{\ell}; p)} \right]
=\frac{E\left( \alpha \prod_{j=1}^N s_j/r_j ; p\right)}
{E(\alpha; p)}
\prod_{1 \leq \ell < m \leq N}
\frac{s_m}{r_m} \frac{E(s_{\ell}/s_m; p)}{E(r_{\ell}/r_m; p)}.
$$
If we use (\ref{eqn:TVa1}), 
we obtain (\ref{eqn:equality1}). \qed
\vskip 0.3cm

The following equality is given
in the first line of Proposition 6.1 in \cite{RS06}.

\begin{lem}
\label{thm:equality2}
For $\s \in \C^N$, $\kappa \in \C$, 
\begin{eqnarray}
&& E \left( \kappa \prod_{j=1}^N s_j ; p \right)
W_{A_{N-1}}(\s; p)
= \frac{(p^N; p^N)_{\infty}^N}{(p;p)_{\infty}^N}
\det_{1 \leq j,k \leq N}
\Big[ s_j^{k-1} E \Big(
(-1)^{N-1} p^{k-1} \kappa s_j^{N}; p^N \Big) \Big].
\nonumber
\end{eqnarray}
\end{lem}

From now on, we assume
$
\Re \tau = 0, \Im \tau > 0$,
that is, $0 < q=e^{\pi i \tau} < 1$.
It is obvious from (\ref{eqn:Jacobi3}) that
if $v \in \R$, then
$\vartheta_1(v; \tau) \in \R, 
|\vartheta_1(v; \tau)| < \infty$.
We set
$p=q^2$
and 
$s_j = e^{i x_j/r}$,
$r_j = e^{i u_j/r}$, $\kappa =e^{i \delta/r}$
in Lemma \ref{thm:equality1}, 
$x_j, u_j \in \R, 1 \leq j \leq N$,
$\delta \in \pi r \Z$.
Then, through (\ref{eqn:JacobiB1}), these lemmas 
are rewritten as follows.
\begin{lem}
\label{thm:equalityB1}
Assume that
$\u=(u_1, \dots, u_N) \in \cA_{[0, 2 \pi r)^N}$ and
$\overline{u}_{\delta} = \delta+\sum_{j=1}^N u_j \in (0, 2 \pi r)$.
Let $\overline{x}_{\delta}=\delta+\sum_{j=1}^N x_j$. 
Then
\begin{eqnarray}
&& \det_{1 \leq j, k \leq N}
\left[ \frac{\vartheta_1((\overline{u}_{\delta}+x_j-u_k)/2 \pi r ; \tau)}
{\vartheta_1(\overline{u}_{\delta}/2 \pi r; \tau)}
\prod_{1 \leq \ell \leq N, \ell \not= k}
\frac{\vartheta_1((x_j-u_{\ell})/2 \pi r; \tau)}
{\vartheta_1((u_k-u_{\ell})/2 \pi r;\tau)} 
\right]
\nonumber\\
&& \quad
= \frac{\vartheta_1(\overline{x}_{\delta}/2 \pi r; \tau)}
{\vartheta_1(\overline{u}_{\delta}/2 \pi r; \tau)}
\prod_{1 \leq j < k \leq N}
\frac{\vartheta_1((x_j-x_k)/2 \pi r; \tau)}
{\vartheta_1((u_j-u_k)/2 \pi r; \tau)}.
\nonumber
\end{eqnarray}
\end{lem}
Similarly 
Lemma \ref{thm:equality2} gives the following.
\begin{lem}
\label{thm:equalityB2}
\begin{eqnarray}
&& \vartheta_1 \left( \frac{\overline{x}_{\delta}}{2 \pi r}; \tau \right)
\prod_{1 \leq j < k \leq N}
\vartheta_1 \left( \frac{x_j-x_k}{2 \pi r}; \tau \right)
\nonumber\\
\label{eqn:equalityB2}
&& \qquad
= C_N^A(\tau)
\det_{1 \leq j, k \leq N}
\left[ e^{i (k-1) x_j/r}
\vartheta_1 \left(
\frac{N-1}{2} +(k-1)\tau+ \frac{\delta+N x_j}{2 \pi r}; N \tau \right) \right]
\end{eqnarray}
with
\begin{eqnarray}
C_N^A(\tau)
&=& C_N^A(\tau; r, \delta)
= q_0(\tau)^{(N-1)(N-2)/2} i^{(N-1)(3N-2)/2} q(\tau)^{(N-1)(3N-2)/8}e^{i (N-1) \delta/2r}
\nonumber\\
\label{eqn:equalityB2b}
&=& q_0(\tau)^{(N-1)(N-2)/2}
\exp \left[
(N-1) \left\{ \frac{(3N-2)}{8} \tau
+\frac{\delta}{2 \pi r}+\frac{3N-2}{4} \right\} \pi i \right].
\end{eqnarray}
\end{lem}

Let $\eta(x)$ denotes Dedekind's $\eta$-function
\cite{Dys72,Mac72}, 
\begin{equation}
\eta(x)=x^{1/24} \prod_{n=1}^{\infty} (1-x^n).
\label{eqn:Dedekind}
\end{equation}
The following equalities were proved as Proposition 5.6.3
in \cite{For10} (see also \cite{For90a,For90b}).
\begin{lem}
\label{thm:Forrester}
Let $\alpha \in \C$.
For $N$ odd
\begin{eqnarray}
&& \det_{1 \leq j, k \leq N} [ \vartheta_3(x_j+\alpha-k/N; \tau)]
= N^{N/2} \eta(e^{2 N \pi i \tau})^{-(N-1)(N-2)/2}
\nonumber\\
&& \qquad \qquad \times
\vartheta_3 \left( \sum_{j=1}^N (x_j+\alpha)+\frac{N \tau}{2}; 2 N \tau \right)
\prod_{1 \leq j < k \leq N} \vartheta_1(x_k-x_j; N \tau),
\label{eqn:Forrester1}
\end{eqnarray}
while for $N$ even
\begin{eqnarray}
&& \det_{1 \leq j, k \leq N} [ \vartheta_1(x_j+\alpha-k/N; \tau)]
= N^{N/2} \eta(e^{2 N \pi i \tau})^{-(N-1)(N-2)/2}
\nonumber\\
&& \qquad \qquad \times
\vartheta_0 \left( \sum_{j=1}^N (x_j+\alpha)+\frac{N \tau}{2}; 2 N \tau \right)
\prod_{1 \leq j < k \leq N} \vartheta_1(x_k-x_j; N \tau).
\label{eqn:Forrester2}
\end{eqnarray}
\end{lem}

\subsection{Generalized $h$-transform \label{sec:h_transform}}

The backward Kolmogorov equation for (\ref{eqn:SDEe1}) is 
given as
\begin{eqnarray}
-\frac{\partial u^A(s,\x)}{\partial s} &=& 
\frac{1}{2} \sum_{j=1}^N \frac{\partial^2 u^A(s,\x)}{\partial x_j^2} 
+ \sum_{\substack{1 \leq j, k \leq N, \cr j \not=k}} A_N^{2 \pi r}(t_{*}-s, x_j-x_k)
\frac{\partial u^A(s,\x)}{\partial x_j}
\nonumber\\
&& + \sum_{1 \leq j \leq N}  A_N^{2 \pi r} (t_{*}-s, \overline{x}_{\delta} )
\frac{\partial u^A(s,\x)}{\partial x_j}.
\label{eqn:Keq1}
\end{eqnarray}
We write the tpd of the process $\Xi^A(t), t \in [0, t_*)$ as
$p_N^A(t, \y|s, \x)=p_N^A(t, \y|s, \x; r, t_*), 0 \leq s \leq t < t_{*}$,
provided that $\x, \y \in \cA_{[0, 2 \pi r)^N}$ 
and $\overline{x}_{\delta}, \overline{y}_{\delta} \in (0, 2 \pi r)$.
Since the configuration of the process $\Xi^A(t), t \in [0, t_*)$
is unlabeled as (\ref{eqn:Xi1}), 
we solve (\ref{eqn:Keq1}) to obtain $p_N^A$ under 
the `initial condition' 
\begin{equation}
\lim_{s \uparrow t} u^A(s, \x)=
\sum_{\sigma \in \cS_N} \prod_{j=1}^N \delta_{y_{\sigma(j)}}( \{x_j \}),
\label{eqn:Keq2}
\end{equation}
where $\cS_N$ denotes a collection of all permutations of $N$ indices.
It is a moderated version of the usual one
$\lim_{s \uparrow t} u^A(s, \x)= \prod_{j=1}^N \delta_{y_j}( \{x_j\})$
for processes with labeled configurations.

We consider the Brownian motion $\V^r(\cdot)$ 
started at $\u \in \cA_{[0, 2 \pi r)^N}$
with an index $\delta \in \pi r \Z$ chosen as
$\overline{u}_{\delta} \in (0, 2 \pi r)$, 
which is killed when it arrives at the boundary of $\cA^{A_{N-1}}_{2 \pi r}$
and when $\overline{V}^{r}_{\delta} (\cdot)\in \{0, 2 \pi r\}$.
Let $q_N^{A}(t-s, \y|\x)=q_N^{A}(t-s, \y|\x; r)$,
$\x, \y \in \cA_{[0, 2 \pi r)^N}$, 
$0 < s < t < t_{*}$ be the tpd of $\V^r(\cdot)$,
which satisfies
\begin{equation}
\lim_{t \downarrow 0} q^A_N(t, \y|\x)= 
\sum_{\sigma \in \cS_N} \prod_{j=1}^N \delta_{y_{\sigma(j)}}( \{x_j \}),
\label{eqn:q_initial}
\end{equation}

\begin{lem}
\label{thm:h_trans}
The tpd of the process $\Xi^A(t), t \in [0, t_*)$ is given by
$$
p_N^A(t, \y|s, \x)
=\frac{h_N^{A}(t_{*}-t, \y)}{h_N^{A}(t_{*}-s, \x)}
q_N^A(t-s, \y|\x),
\quad 0 \leq s \leq t < t_{*},  \quad \x, \y \in \cA_{[0, 2 \pi r)^N}, 
$$
where $\overline{x}_{\delta}, \overline{y}_{\delta} \in (0, 2 \pi r)$ and
\begin{eqnarray}
h_N^A(t_{*}-t,\x) &=& h_N^A(t_{*}-t,\x; r, t_*)
\nonumber\\
&=& 
e^{-N(N-1)(N-2) t_{*}/48 r^2}
\eta(e^{-N (t_{*}-t)/r^2})^{-(N-1)(N-2)/2}
\nonumber\\
&& \times 
\vartheta_1 \left(\frac{\overline{x}_{\delta}}{2 \pi r} ; 
\frac{iN(t_{*}-t)}{2 \pi r^2} \right)
\prod_{1 \leq j < k \leq N}
\vartheta_1 \left(\frac{x_k-x_j}{2 \pi r} ; \frac{iN(t_{*}-t)}{2\pi r^2} \right),
\label{eqn:h1}
\end{eqnarray}
$t \in [0, t_{*}), \x \in \cA_{[0, 2 \pi r)^N}$.
\end{lem}
\vskip 0.3cm
\noindent{\it Proof.} 
Set
\begin{equation}
u^A(s,\x)=f^A(s, \x) q_N^{A}(t-s,\y|\x)
\label{eqn:eqB1}
\end{equation}
and put it into (\ref{eqn:Keq1})
assuming that $f^A$ is $\rC^1$ in $t$ and $\rC^2$ in $\x$.
Then we have
\begin{eqnarray}
&& -\frac{\partial f^A(s,\x)}{\partial s} q_N^{A}(t-s, \y|\x)
\nonumber\\
&& = \frac{1}{2} q_N^{A}(t-s,\y|\x) \sum_{j=1}^N \frac{\partial^2 f^A(s,\x)}{\partial x_j^2}
+\sum_{j=1}^N \frac{\partial f^A(s,\x)}{\partial x_j}
\frac{\partial q_N^{A}(t-s,\y|\x)}{\partial x_j}
\nonumber\\
&& + q_N^{A}(t-s,\y|\x) 
\sum_{j=1}^N \left\{ \sum_{\substack{1 \leq k \leq N, \cr k \not= j}}
A_N^{2 \pi r}(t_{*}-s, x_j-x_k) 
+ A_N^{2 \pi r} (t_{*}-s, \overline{x}_{\delta} ) \right\}
\frac{\partial f^A(s,\x)}{\partial x_j}
\nonumber\\
&& + f^A(s,\x) 
\sum_{j=1}^N \left\{ \sum_{\substack{1 \leq k \leq N, \cr k \not= j}}
A_N^{2 \pi r}(t_{*}-s, x_j-x_k) 
+  A_N^{2 \pi r} (t_{*}-s, \overline{x}_{\delta}) \right\}
\frac{\partial q_N^{A}(t-s,\y|\x)}{\partial x_j},
\label{eqn:eqB2}
\end{eqnarray}
since $q_N^{A}(t-s,\y|\x)$ satisfies the diffusion equation.
We put
\begin{equation}
f^A(s,\x)= g^A(s)
\left\{ \vartheta_1 \left(\frac{\overline{x}_{\delta}}{2 \pi r} ; 
\frac{iN(t_{*}-s)}{2 \pi r^2} \right)
\prod_{1 \leq j < k \leq N}
\vartheta_1 \left(\frac{x_k-x_j}{2 \pi r} ; \frac{iN(t_{*}-s)}{2\pi r^2} \right) \right\}^{-1}, 
\label{eqn:f1}
\end{equation}
where $g^A$ is a $\rC^1$ function of time $s$ to be determined. 
By definition (\ref{eqn:A2}), we see
$$
\frac{\partial f^A(s,\x)}{\partial x_j}
= - \left\{ \sum_{\substack{1 \leq k \leq N, \cr k \not=j}}
A_N^{2 \pi r}(t_{*}-s, x_j-x_k)
+  A_N^{2 \pi r} (t_{*}-s, \overline{x}_{\delta} ) 
\right\}f^A(s,\x)
$$
and
\begin{eqnarray}
&& \frac{\partial A_N^{2 \pi r}(t_{*}-s, x_j-x_k)}{\partial x_j}
\nonumber\\
&& \quad = \frac{1}{(2 \pi r)^2}
\frac{\vartheta_1''((x_j-x_k)/2\pi r; iN(t_{*}-s)/2\pi r^2)}
{\vartheta_1((x_j-x_k)/2 \pi r; iN(t_{*}-s)/2\pi r^2)}
-(A_N^{2 \pi r}(t_{*}-s, x_j-x_k))^2.
\label{eqn:f5}
\end{eqnarray}
Then (\ref{eqn:eqB2}) gives the equation
\begin{eqnarray}
-\frac{\partial f^A(s,\x)}{\partial s}
&=& - \frac{1}{2(2 \pi r)^2}
\sum_{\substack{1 \leq j, k \leq N, \cr j \not=k}}
\frac{\vartheta_1''((x_j-x_k)/2 \pi r; iN(t_{*}-s)/2\pi r^2)}
{\vartheta_1((x_j-x_k)/2 \pi r; iN(t_{*}-s)/2\pi r^2)} f^A(s,\x)
\nonumber\\
&& - \frac{1}{2(2 \pi r)^2} N
\frac{\vartheta_1''(\overline{x}_{\delta}/2 \pi r; iN(t_{*}-s)/2\pi r^2)}
{\vartheta_1(\overline{x}_{\delta}/2 \pi r; iN(t_{*}-s)/2\pi r^2)} f^A(s,\x)
\nonumber\\
&& - \frac{1}{2} 
\sum_{\substack{1 \leq j, k, \ell \leq N, \cr j \not= k \not= \ell}}
A_N^{2 \pi r}(t_{*}-s, x_j-x_k) A_N^{2 \pi r}(t_{*}-s, x_j-x_{\ell}) f^A(s,\x),
\label{eqn:eqB5}
\end{eqnarray}
where the sum in the last term denotes
the summation over $1 \leq j, k, \ell \leq N$ conditioned that
$j, k, \ell$ are all distinct.
By the setting (\ref{eqn:f1}), 
\begin{eqnarray}
\mbox{LHS of (\ref{eqn:eqB5})}
&=& -\frac{1}{g^A(s)} \frac{d g^A(s)}{ds} f^A(s,\x)
\nonumber\\
&& - \frac{iN}{4 \pi r^2} 
\sum_{\substack{1 \leq j, k \leq N, \cr j \not= k}}
\frac{\dot{\vartheta_1}((x_j-x_k)/2 \pi r; iN(t_{*}-s)/2\pi r^2)}
{\vartheta_1((x_j-x_k)/2 \pi r; iN(t_{*}-s)/2\pi r^2)} f^A(s,\x)
\nonumber\\
&& - \frac{iN}{2\pi r^2} 
\frac{\dot{\vartheta_1}(\overline{x}_{\delta}/2 \pi r; iN(t_{*}-s)/2\pi r^2)}
{\vartheta_1(\overline{x}_{\delta}/2 \pi r; iN(t_{*}-s)/2\pi r^2)} f^A(s,\x), 
\nonumber
\end{eqnarray}
where $\dot{\vartheta}_1(x; \tau)=d \vartheta_1(x; \tau)/d \tau$
and (\ref{eqn:Jacobi4}) was used.
Since (\ref{eqn:Jacobi_eq}) is satisfied, the above is equal to
\begin{eqnarray}
-\frac{1}{g^A(s)} \frac{d g^A(s)}{ds} f^A(s,\x)
&-& \frac{N}{16 \pi^2 r^2} 
\sum_{\substack{1 \leq j, k \leq N, \cr j \not= k}}
\frac{\vartheta_1''((x_j-x_k)/2 \pi r; iN(t_{*}-s)/2\pi r^2)}
{\vartheta_1((x_j-x_k)/2 \pi r; iN(t_{*}-s)/2\pi r^2)} f^A(s,\x)
\nonumber\\
&-& \frac{N}{8 \pi^2 r^2} 
\frac{\vartheta_1''(\overline{x}_{\delta}/2 \pi r; iN(t_{*}-s)/2\pi r^2)}
{\vartheta_1(\overline{x}_{\delta}/2 \pi r; iN(t_{*}-s)/2\pi r^2)} f^A(s,\x).
\nonumber
\end{eqnarray}
Therefore, (\ref{eqn:eqB5}) becomes
\begin{eqnarray}
\frac{1}{g^A(s)} \frac{d g^A(s)}{ds}
&=& \frac{1}{2}
\sum_{\substack{1 \leq j, k, \ell \leq N, \cr j \not= k \not= \ell}}
A_N^{2 \pi r}(t_{*}-s, x_j-x_k) A_N^{2 \pi r}(t_{*}-s, x_j-x_{\ell})
\nonumber\\
&& - \frac{N-2}{16 \pi^2 r^2}
\sum_{\substack{1 \leq j, k \leq N, \cr j \not= k}}
\frac{\vartheta_1''((x_j-x_k)/2 \pi r; iN(t_{*}-s)/2\pi r^2)}
{\vartheta_1((x_j-x_k)/2 \pi r; iN(t_{*}-s)/2\pi r^2)}.
\label{eqn:eqB6}
\end{eqnarray}
Now we rewrite RHS of (\ref{eqn:eqB6}) 
by using the functions $\wp$ and $\zeta$
through (\ref{eqn:A1}).
First we see
\begin{eqnarray}
&& \frac{1}{2} \sum_{\substack{1 \leq j, k, \ell \leq N, \cr j \not= k \not= \ell}}
A_N^{2 \pi r}(t_{*}-s, x_j-x_k) A_N^{2 \pi r}(t_{*}-s, x_j-x_{\ell})
\nonumber\\
&& \qquad = \frac{1}{2} \sum_{\substack{1 \leq j, k, \ell \leq N, \cr j \not= k \not= \ell}}
\zeta(x_j-x_k)\zeta(x_j-x_{\ell})
-\frac{\eta_1(t_{*}-s)}{\pi r} 
\sum_{\substack{1 \leq j, k, \ell \leq N, \cr j \not= k \not= \ell}}
\zeta(x_j-x_k)(x_j-x_{\ell})
\nonumber\\
&& \qquad \qquad
+ \frac{1}{2} \left( \frac{\eta_1(t_{*}-s)}{\pi r} \right)^2 
\sum_{\substack{1 \leq j, k, \ell \leq N, \cr j \not= k \not= \ell}}
(x_j-x_k)(x_j-x_{\ell}).
\nonumber
\end{eqnarray}
For
\begin{eqnarray}
&&
 \sum_{\substack{1 \leq j, k, \ell \leq N, \cr j \not= k \not= \ell}}
\zeta(x_j-x_k)\zeta(x_j-x_{\ell})
\nonumber\\
&& \quad 
=2 \sum_{1 \leq j < k < \ell \leq N}
\Big\{ \zeta(x_j-x_k)\zeta(x_j-x_{\ell})
+\zeta(x_k-x_{\ell})\zeta(x_k-x_j)
+\zeta(x_{\ell}-x_j)\zeta(x_{\ell}-x_k) \Big\},
\nonumber
\end{eqnarray}
Lemma \ref{thm:zeta_P1} gives
\begin{eqnarray}
&&
 \frac{1}{2} \sum_{\substack{1 \leq j, k, \ell \leq N, \cr j \not= k \not= \ell}}
\zeta(x_j-x_k)\zeta(x_j-x_{\ell})
=  \frac{1}{4} 
\sum_{\substack{1 \leq j, k, \ell \leq N, \cr j \not= k \not= \ell}}
\zeta(x_j-x_k)^2
-\frac{1}{4}
\sum_{\substack{1 \leq j, k, \ell \leq N, \cr j \not= k \not= \ell}}
\wp(x_j-x_k)
\nonumber\\
&& \qquad \qquad \qquad 
= \frac{N-2}{4} 
\sum_{\substack{1 \leq j, k \leq N, \cr j \not= k}}
\zeta(x_j-x_k)^2
-\frac{N-2}{4}
\sum_{\substack{1 \leq j, k \leq N, \cr j \not= k}}
\wp(x_j-x_k).
\nonumber
\end{eqnarray}
On the other hand, differentiation of (\ref{eqn:A1}) with respect to $x$
gives
$$
\frac{\partial A_N^{2 \pi r}(t_{*}-s,x)}{\partial x}
= \zeta'(x)-\frac{\eta_1(t_{*}-s)}{\pi r}
= - \wp(x)-\frac{\eta_1(t_{*}-s)}{\pi r},
$$
where (\ref{eqn:wp2}) was used.
Combining it with (\ref{eqn:f5}) gives
\begin{eqnarray}
&& \frac{1}{(2 \pi r)^2}
\frac{\vartheta_1''(x/2 \pi r; iN(t_*-s)/2 \pi r^2)}{\vartheta_1(x/2 \pi r; iN(t_*-s)/2 \pi r^2)}
= A_N^{2 \pi r}(t_{*}-s,x)^2-\wp(x)-\frac{\eta_1(t_{*}-s)}{\pi r}
\nonumber\\
&& \qquad \qquad 
= \zeta(x)^2- \frac{2 \eta_1(t_{*}-s)}{\pi r} \zeta(x) x 
+\left( \frac{\eta_1(t_{*}-s)}{\pi r} \right)^2 x^2
-\wp(x)-\frac{\eta_1(t_{*}-s)}{\pi r}.
\nonumber
\end{eqnarray}
Then
\begin{eqnarray}
&& - \frac{N-2}{16 \pi^2 r^2}
\sum_{\substack{1 \leq j, k \leq N, \cr j \not= k}}
\frac{\vartheta_1''((x_j-x_k)/2 \pi r; iN(t_{*}-s)/2\pi r^2)}
{\vartheta_1((x_j-x_k)/2 \pi r; iN(t_{*}-s)/2\pi r^2)}
\nonumber\\
&& \quad =
- \frac{N-2}{4} 
\sum_{\substack{1 \leq j, k \leq N, \cr j \not= k}}
\zeta(x_j-x_k)^2
+\frac{N-2}{2} \frac{\eta_1(t_{*}-s)}{\pi r} 
\sum_{\substack{1 \leq j, k \leq N, \cr j \not= k}}
\zeta(x_j-x_k) (x_j-x_k)
\nonumber\\
&& \quad \quad 
- \frac{N-2}{4} \left( \frac{\eta_1(t_{*}-s)}{\pi r} \right)^2
\sum_{\substack{1 \leq j, k \leq N, \cr j \not= k}}
(x_j-x_k)^2
+\frac{N-2}{4} \sum_{\substack{1 \leq j, k \leq N, \cr j \not= k}}
\wp(x_j-x_k)
\nonumber\\
&& \qquad 
+ \frac{\eta_1(t_{*}-s)}{\pi r} \frac{N(N-1)(N-2)}{4}.
\nonumber
\end{eqnarray}
It is easy to prove that
\begin{eqnarray}
&&
\sum_{\substack{1 \leq j, k, \ell \leq N, \cr j \not= k \not= \ell}}
\zeta(x_j-x_k)(x_j-x_{\ell})
-\frac{N-2}{2}
\sum_{\substack{1 \leq j, k \leq N, \cr j \not= k}}
\zeta(x_j-x_k)(x_j-x_k)
= 0,
\nonumber\\
&& 
\sum_{\substack{1 \leq j, k, \ell \leq N, \cr j \not= k \not= \ell}}
(x_j-x_k)(x_j-x_{\ell})
-\frac{N-2}{2}
\sum_{\substack{1 \leq j, k \leq N, \cr j \not= k}}
(x_j-x_k)^2
= 0, 
\nonumber
\end{eqnarray}
by using the fact that $\zeta$ is odd.
Then the equation (\ref{eqn:eqB6}) is reduced to be
$$
\frac{d}{ds} \log g^A(s)
=\frac{\eta_1(t_{*}-s)}{4 \pi r} N(N-1)(N-2).
$$
Since $\eta_1(t_{*}-t)$ is explicitly given as (\ref{eqn:eta1b}) with $\alpha=2 \pi r$,
this equation can be solved as
$$
g^A(s)=c'e^{N(N-1)(N-2) s/48r^2}
 \prod_{n=1}^{\infty}
\left( \frac{1-e^{-nN(t_{*}-s)/r^2}}{1-e^{-nN t_{*}/r^2}}
\right)^{(N-1)(N-2)/2}
$$
with a constant $c'$.
The solution (\ref{eqn:eqB1}) has been determined of the form
$$
u^A(s,\x)
= \prod_{n=1}^{\infty}
(1-e^{-nN t_{*}/r^2})^{-(N-1)(N-2)/2}
\frac{c'}{h_N^{A}(t_{*}-s,\x)}
q_N^{A}(t-s, \y|\x).
$$
For (\ref{eqn:q_initial}), 
the condition (\ref{eqn:Keq2}) is satisfied, if and only if
$$
c'=
\prod_{n=1}^{\infty}
(1-e^{-nN t_{*}/r^2})^{(N-1)(N-2)/2}
h_N^{A}(t_{*}-t,\y).
$$
Since $\vartheta_1(x/2\pi r; iN(t_{*}-t)/2 \pi r^2)>0$
if $x \in (0, 2 \pi r)$, $t \in [0, t_{*})$, 
and $q_N^A(t-s, \y|\x)$ is assumed to be the tpd of $\V^r(\cdot)$, 
we can conclude that 
$0< p_N^{A}(t, \y|s, \x) < \infty$ for any
$\x, \y \in \cA_{[0, 2 \pi r)^N}$, 
$0 \leq s \leq t < t_{*}$, 
where $\delta$ is chosen so that
$\overline{x}_{\delta}, \overline{y}_{\delta} \in (0, 2 \pi r)$.
Then the proof is completed. \qed
\vskip 0.3cm

Let $\check{\bW}(t)=(\chW_1(t), \dots, \chW_N), t \geq 0$
be $N$-dimensional Brownian motion on $(\rS^1(r))^N$
started at $\u \in \cA^{A_{N-1}}_{2 \pi r}$.
The expectation with respect to this process is denoted 
by $\check{\rE}_{\u}$.
Consider a stopping time
$$
T_{\check{\bbW}}=\inf\{t > 0: \check{\bW}(t) \notin \cA^{A_{N-1}}_{2 \pi r}\}.
$$
Put $\overline{W}_{\delta}=\delta+\sum_{j=1}^N \chW_j(t)$,
where the index $\delta \in \pi r \Z$ is determined
so that $\overline{u}_{\delta} \in (0, 2 \pi r)$.
Then we also consider the following stopping time
$$
T_{\overline{W}_{\delta}}=\inf\{t > 0: \overline{W}_{\delta} \in \{0, 2 \pi r\} \}.
$$

For the process $\Xi^A(t), t \in [0, t_{*})$ is
$\mM([0, 2 \pi r))$-valued,
measurable functions are symmetric functions of 
$N$ variables $X^A_j, 1 \leq j \leq N$ at each time.
By the definition (\ref{eqn:Def_Xjs}) for $\X^A$,
they should be periodic with period $2 \pi r$.
Let $T \in [0, t_{*})$.
Then any $\cF_{\Xi^A}(T)$-measurable function $F$
will be given as follows.
With an arbitrary integer $M \in \N$ and
arbitrary sequence of times $0 \leq t_1 < \cdots < t_M \leq T$,
\begin{equation}
F(\Xi^A(\cdot))=\prod_{m=1}^M g_m(\X^A(t_m)),
\label{eqn:measure1}
\end{equation}
where $g_m(\x), 1 \leq m \leq M$ are symmetric functions
and
\begin{equation}
g_m((x_j+2 \pi r n_j))=g_m(\x), \quad
n_j \in \Z, \quad 1 \leq j \leq N,
\label{eqn:g_periodic}
\end{equation}
for $1 \leq m \leq M$.

The indicator function of $\omega$ is
denoted by $\1(\omega)$; 
$\1(\omega)=1$ if $\omega$ is satisfied, and
$\1(\omega)=0$ otherwise.
Lemma \ref{thm:h_trans} implies the following equality.

\begin{prop}
\label{thm:h_trans2}
Suppose $\xi=\sum_{j=1}^N \delta_{u_j} \in \mM_0([0, 2 \pi r))$.
Let $T \in [0, t_{*})$.
For any $\cF_{\Xi^A}(T)$-measurable observable $F$,
$$
\E_{\xi}^A[F(\Xi^A(\cdot))]
=\check{\rE}_{\u} \left[ F \left( \sum_{j=1}^N \delta_{\chW_j(\cdot)} \right)
\1(T_{\check{\bbW}} \wedge T_{\overline{W}_{\delta}} > T)
\frac{h_N^{A}(t_{*}-T, \check{\bW}(T))}{h_N^{A}(t_{*}, \u)} \right].
$$
\end{prop}
\vskip 0.3cm
See Remark 2 at the end of Section \ref{sec:DMR}.

\subsection{Markov process $\bW^r$  \label{sec:bWr}}

We write the tpd of BM on $\R$ as
$$
p_{\rm BM}(t, y|x)=\frac{1}{\sqrt{2 \pi t}}
e^{-(y-x)^2/2t},
\quad x, y \in \R, \quad t \in [0, \infty).
$$
By wrapping it on $\rS^1(r)$, we define
\begin{equation}
p^r_{A_{N-1}}(t,y|x) = \left\{
\begin{array}{ll}
\displaystyle{
\sum_{\ell \in \Z} p_{\rm BM}(t, y+2 \pi r \ell|x)
}, 
 & \mbox{if $N$ is even}, \cr
\displaystyle{
\sum_{\ell \in \Z} (-1)^{\ell}
p_{\rm BM}(t, y+2 \pi r \ell|x)
}, 
\quad & \mbox{if $N$ is odd}, 
\end{array} \right.
\label{eqn:prZ1}
\end{equation}
$x, y \in [0, 2 \pi r)$, $t \geq 0$.
Using the Jacobi theta functions (\ref{eqn:theta_0_2_3}), 
it is written as
$$
p^r_{A_{N-1}}(t,y|x) = \left\{
\begin{array}{ll}
\displaystyle{
p_{\rm BM}(t, y|x) \vartheta_3 \left( \frac{i(y-x)r}{t}; \frac{2 \pi i r^2}{t} \right)
}, \quad & \mbox{if $N$ is even}, \cr
\displaystyle{
p_{\rm BM}(t, y|x) \vartheta_0 \left( \frac{i(y-x)r}{t}; \frac{2 \pi i r^2}{t} \right)
}, \quad & \mbox{if $N$ is odd}.
\end{array} \right.
$$
We find that by Jacobi's imaginary transformations \cite{WW27},
\begin{eqnarray}
\vartheta_0(v; \tau)
&=& e^{\pi i/4} \tau^{-1/2} e^{-\pi i v^2/\tau}
\vartheta_2 \left( \frac{v}{\tau}; - \frac{1}{\tau} \right),
\nonumber\\
\vartheta_3(v; \tau)
&=& e^{\pi i/4} \tau^{-1/2} e^{-\pi i v^2/\tau}
\vartheta_3 \left( \frac{v}{\tau}; - \frac{1}{\tau} \right),
\nonumber
\end{eqnarray}
the above is further rewritten as
\begin{equation}
p^r_{A_{N-1}}(t,y|x) = \left\{
\begin{array}{ll}
\displaystyle{
\frac{1}{2 \pi r} \vartheta_3 \left( \frac{y-x}{2 \pi r}; \frac{it}{2 \pi r^2} \right)
}, \quad & \mbox{if $N$ is even}, \cr
\displaystyle{
\frac{1}{2 \pi r} \vartheta_2 \left( \frac{y-x}{2 \pi r}; \frac{it}{2 \pi r^2} \right)
}, \quad & \mbox{if $N$ is odd}.
\end{array} \right.
\label{eqn:prZ3}
\end{equation}

Lemma \ref{thm:Forrester} given by Forrester \cite{For10}
implies the following.

\begin{prop}
\label{thm:qNA}
For $N \in \{2,3, \dots \}$, $\v$ is  given by (\ref{eqn:ed1}).
Then for $\y \in \cA_{[0, 2 \pi r)^N}$,
$t>0$,
\begin{equation}
q_N^A(t, \y|\v)
=\det_{1 \leq j, k \leq N} 
\Big[ p^r_{A_{N-1}}(t, y_j|v_k) \Big].
\label{eqn:qNA_1}
\end{equation}
\end{prop}
\noindent{\it Proof.} 
For $N$ odd, we put $\alpha=1/N+\tau/2$ in (\ref{eqn:Forrester1}),
and for $N$ even, $\alpha=1/N+(1+\tau)/2$ in (\ref{eqn:Forrester2}).
Let $\tau=it/2 \pi r^2$. Then we have
\begin{eqnarray}
&& \det_{1 \leq j, k \leq N}
\Big[ p^r_{A_{N-1}}(t, y_j| v_k) \Big] =
\left( \frac{\sqrt{N}}{2 \pi r} \right)^N
\eta(e^{-Nt/r^2})^{-(N-1)(N-2)/2} e^{Nt/8r^2}
\nonumber\\
&& \qquad \qquad \times
\vartheta_1 \left( \frac{\overline{y}_{- \pi r (N-2)}}{2\pi r} ;
\frac{i N t}{\pi r^2} \right)
\prod_{1 \leq j < k \leq N}
\vartheta_1 \left( \frac{y_k-y_j}{2 \pi r} ; \frac{i N t}{2 \pi r^2} \right),
\label{eqn:qNA_2}
\end{eqnarray}
where quasi-periodicity of $\vartheta_{\mu}, \mu=0,1,2,3$
has been used.
By expression (\ref{eqn:prZ3}) with (\ref{eqn:Jacobi_equation}),
it is obvious that (\ref{eqn:qNA_2}) satisfies the diffusion equation.
This expression (\ref{eqn:qNA_2}) guarantees the positivity and finiteness
of $\det_{1 \leq j, k \leq N}[ p^r_{A_{N-1}}(t, y_j| v_k)]$
for $\y \in \cA_{[0, 2 \pi r)^N}$ and
$\overline{y}_{-\pi r (N-2)} \in (0, 2 \pi r)$.
Equation (\ref{eqn:qNA_2}) also shows that it vanishes when
$y_j=y_k$ for any $j \not= k$ and when
$\overline{y}_{-\pi r (N-2)} \in \{0, 2\pi r\}$.
By the expression (\ref{eqn:prZ1}) 
and the argument given by Liechty and Wang \cite{LW13}
(see also \cite{Ful04}), 
we can prove that
(\ref{eqn:qNA_2}) satisfies the moderated initial configuration
$$
\lim_{t \downarrow 0}
\det_{1 \leq j, k \leq N}
[ p^r_{A_{N-1}}(t, y_j|v_k)]
=\sum_{\sigma \in \cS_N} \prod_{j=1}^N \delta_{v_{\sigma(j)}}(\{y_j\}).
$$
Then the proof is completed. \qed
\vskip 0.3cm
\noindent{\bf Remark 1.} \,
The Karlin-McGregor-type \cite{KM59}
determinantal formula (\ref{eqn:qNA_1}) for $q_N^A(t, \y|\x)$
given in Proposition \ref{thm:qNA} is crucial for DMR which will be
proved in Theorem \ref{thm:DMR}.
In the present paper, we obtained it for the special 
initial configuration (\ref{eqn:ed1})
due to the explicit evaluation (\ref{eqn:qNA_2}) given by
Lemma \ref{thm:Forrester} of Forrester \cite{For10}.
If we obtain the Karlin-McGregor-type determinantal formula
for $q_N^A$ for other initial configuration,
we can prove DMR for the process and it will immediately
conclude that the process is determinantal
by Theorem 1.3 of \cite{Kat13a}.
\vskip 0.3cm

Given $N \in \N$, let
$W^r(t), t \geq 0$ be a Markov process in $[0, 2 \pi r)$
such that its transition density is given by
$p^r_{A_{N-1}}(t, y|x), t \geq 0, x, y \in [0, 2 \pi r)$,
defined by (\ref{eqn:prZ1}).
Then we introduce an $N$ independent copies of $W^r(t), t \geq 0$,
denoted by $W^r_j(t), t \geq 0, 1 \leq j \leq N$ and
let $\bW^r(t)=(W^r_1(t), \dots, W^r_N(t))$, $t \geq 0$.
The probability space of the process is denoted by
$(\Omega_{W^r}, \cF_{W^r}, \rP_{\v}^r)$,
and the expectation is written as $\rE_{\v}^r$,
where the initial configuration is given by $\v$ with (\ref{eqn:ed1}). 
A filtration $\{\cF_{W^r}(t) : t \geq 0\}$ is 
generated by $\bW^r(t), t \geq 0$, 
which satisfies the usual conditions. 

\subsection{Martingales and complex Brownian motions \label{sec:CBM}}

Let $0 < t_{*} < \infty$ and
$\xi=\sum_{j=1}^N \delta_{u_j} \in \mM_0([0, 2 \pi r))$.
For $1 \leq k \leq N$, define
\begin{eqnarray}
&& \Phi^A_{\xi, u_k}(z)
= \Phi^A_{\xi, u_k}(z; N, r, t_{*})
\nonumber\\
&& =
\frac{\vartheta_1((\overline{u}_{\delta}+z-u_k)/2 \pi r; iNt_{*}/2 \pi r^2)}
{\vartheta_1(\overline{u}_{\delta}/2 \pi r; iNt_{*}/2 \pi r^2)}
\prod_{\substack{1 \leq \ell \leq N, \cr \ell \not=k}}
\frac{\vartheta_1((z-u_{\ell})/2 \pi r; iNt_{*}/2 \pi r^2)}{\vartheta_1((u_k-u_{\ell})/2 \pi r; iNt_{*}/2 \pi r^2)}, 
\, z \in \C,
\label{eqn:Phi1}
\end{eqnarray}
and
\begin{eqnarray}
\cM^A_{\xi, u_k}(t, x)
&=& \cM^A_{\xi, u_k}(t, x; N, r, t_{*})
\nonumber\\
&=& \int_{\R} d w \,
\frac{e^{-(ix+w)^2/2t}}{\sqrt{2 \pi t}}
\Phi^A_{\xi, u_k}(iw),
\quad (t, x) \in [0, t_{*}) \times [0, 2 \pi r).
\label{eqn:M1}
\end{eqnarray}
Since $\Phi^A_{\xi, u_k}(z), 1 \leq k \leq N$ are
holomorphic for $|z| < \infty$,
(\ref{eqn:M1}) is written as
\begin{eqnarray}
\cM^A_{\xi, u_k}(t,x) &=& \int_{\R} d \widetilde{w} \,
\frac{e^{-\widetilde{w}^2/2t} }{\sqrt{ 2 \pi t}}
\Phi^A_{\xi, u_k}(x+i \widetilde{w})
\nonumber\\
&=& \widetilde{\rE} [ \Phi^A_{\xi, u_k}(x+i \widetilde{W}(t))],
\label{eqn:M2}
\end{eqnarray}
where $\widetilde{W}$ denotes a BM on $\R$ started at 0,
which is independent of $\bW^r$, and
$\widetilde{\rE}$ does the expectation for $\widetilde{W}$.
Then the following is proved.
\begin{lem}
\label{thm:martingale}
Assume $v_j, 1 \leq j \leq N$ are given by (\ref{eqn:ed1})
and $\eta=\sum_{j=1}^N \delta_{v_j}$. 
Then
\begin{eqnarray}
&{\rm (i)}& 
\cM^A_{\eta, v_k}(t, W^r(t)), 1 \leq k \leq N, t \in [0, t_*)
\, \, \mbox{are continuous-time martingales};
\nonumber\\
&& \qquad \qquad 
\rE^r[ \cM^A_{\eta, v_k}(t, W^r(t)) | \cF_{W^r}(s)]
=\cM^A_{\eta, v_k}(s, W^r(s))
\quad \mbox{a.s.}
\nonumber\\
&& \qquad 
\mbox{for any two bounded stopping times
with $0 \leq s \leq t < t_*$.}
\nonumber\\
&{\rm (ii)}& 
\mbox{For any $t \in [0, t_*)$}, 
\cM^A_{\eta, v_k}(t, x),1 \leq k \leq N,
\, \mbox{are linearly independent functions of $x \in [0, 2 \pi r)$},
\nonumber\\
&{\rm (iii)}& 
\cM^A_{\eta, v_k}(0,v_j)=\delta_{jk}, \quad
1 \leq j, k \leq N.
\nonumber
\end{eqnarray}
\end{lem}
\noindent{\it Proof.} 
(i) For the quasi-periodicity (\ref{eqn:Jacobi5a}) of $\vartheta_1$,
the expression (\ref{eqn:M2}) with (\ref{eqn:Phi1}) implies that, 
for $\ell \in \Z$,
\begin{eqnarray}
\cM^A_{\eta, v_k}(t, x+2 \pi r \ell)
&=& (-1)^{\ell N} \cM^A_{\eta, v_k}(t, y)
\nonumber\\
&=& 
\left\{
\begin{array}{ll}
\cM^A_{\eta, v_k}(t, y), \quad & \mbox{if $N$ is even}, \cr
& \cr
(-1)^{\ell} \cM^A_{\eta, v_k}(t, y), \quad & \mbox{if $N$ is odd}.
\end{array} \right.
\label{eqn:Mparity}
\end{eqnarray}
Then, for $0 \leq s \leq t < t_*, 1 \leq k \leq N$,
(\ref{eqn:prZ1}) gives
\begin{eqnarray}
&& \rE^r[\cM^A_{\eta, v_k}(t, W^r(t)) | \cF_{W^r}(s)]
= \int_{0}^{2 \pi r} dw \,
\cM^A_{\eta, v_k}(t, w) p^r_{A_{N-1}}(t-s, w|W^r(s))
\nonumber\\
&& =
\left\{
\begin{array}{ll}
\displaystyle{
\sum_{\ell \in \Z} \int_{2 \pi r \ell}^{2 \pi r (\ell+1)} dw \,
\cM^A_{\eta, v_k}(t, w-2 \pi r \ell) p_{\rm BM}(t-s, w|W^r(s))
}, \quad & \mbox{if $N$ is even}, \cr
\displaystyle{
\sum_{\ell \in \Z} \int_{2 \pi r \ell}^{2 \pi r (\ell+1)} dw \, (-1)^{\ell}
\cM^A_{\eta, v_k}(t, w-2 \pi r \ell) p_{\rm BM}(t-s, w|W^r(s))
}, \quad & \mbox{if $N$ is odd}.
\end{array} \right.
\nonumber
\end{eqnarray}
By (\ref{eqn:Mparity}), it is equal to
$$
\int_{\R} dw \, \cM^A_{\eta, v_k}(t, w) 
p_{\rm BM}(t-s, w|W^r(s)) \quad
\mbox{a.s.}
$$
By definition (\ref{eqn:Jacobi3}) of $\vartheta_1$,
we will obtain the following expansions; for $1 \leq k \leq N$, 
\begin{eqnarray}
\frac{\vartheta_1((\overline{u}_{\delta}+z-u_k)/2 \pi r; iNt_{*}/2 \pi r^2)}
{\vartheta_1(\overline{u}_{\delta}/2 \pi r; iNt_{*}/2 \pi r^2)}
&=& \sum_{n_0 \in \Z} b^0_{n_0} e^{i(2n_0-1) z/2r},
\nonumber\\
\frac{\vartheta_1((z-u_{\ell})/2 \pi r; iNt_{*}/2 \pi r^2)}{\vartheta_1((u_k-u_{\ell})/2 \pi r; iNt_{*}/2 \pi r^2)}
&=& \sum_{n_{\ell} \in \Z} b^{\ell}_{n_{\ell}} e^{i(2n_{\ell}-1) z/2r},
\quad 1 \leq \ell \leq N, \quad \ell \not=k,
\nonumber
\end{eqnarray}
where the coefficients $b^{\ell}_{n_{\ell}}, 0 \leq \ell \leq N, \ell \not= k$
are functions of $\{u_{\ell}\}_{\ell=1}^N$, $t_*, N$ and $r$.
Then, if we introduce an $N$-component index 
$\n=(n_0, n_1, \dots, n_{k-1}, n_{k+1}, \dots, n_N)$ for each $1 \leq k \leq N$,
and put $B^k_{\n}=\prod_{0 \leq \ell \leq N, \ell \not= k} b^{\ell}_{n_{\ell}}$,
(\ref{eqn:M2}) with $\xi=\eta, u_k=v_k, x=w$ is expanded as
\begin{eqnarray}
\cM^A_{\eta, v_k}(t, w)
&=& \sum_{\n \in \Z^N} B^k_{\n} 
\exp \left(i \sum_{0 \leq \ell \leq N, \ell \not=k}
(2n_{\ell}-1) \frac{w}{2r} \right)
\widetilde{\rE} \left[
\exp \left( - \sum_{0 \leq \ell \leq N, \ell \not=k}
(2n_{\ell}-1) \frac{\widetilde{W}(t)}{2r} \right)
\right]
\nonumber\\
&=& \sum_{\n \in \Z^N} B^k_{\n}
G \left( \frac{i}{2r} \sum_{0 \leq \ell \leq N, \ell \not=k} (2n_{\ell}-1) ;
t, w \right),
\label{eqn:M_expantion}
\end{eqnarray}
where
$$
G(\alpha; t, w)=e^{\alpha w-\alpha^2t/2}, \quad \alpha \in \C.
$$
For any $\alpha \in \C$, it is easy to confirm that
$$
\int_{\R} dw \, G(\alpha; t, w) p_{\rm BM}(t-s, w|x)
=G(\alpha; s, x),
\quad 0 \leq s \leq t, \quad x \in \R.
$$
Then (\ref{eqn:M_expantion}) gives
$$
\int_{\R} dw \, \cM^A_{\eta, v_k}(t, w)
p_{\rm BM}(t-s, w|x)
= \cM^A_{\eta, v_k}(s, x),
\quad 0 \leq s \leq t, \quad x \in \R, 
$$
and hence (i) is concluded.
As a matter of course,  $\eta \in \mM_0([0, 2 \pi r))$, and then
the zeroes of $\Phi^A_{\eta, v_j}(z)$ are distinct
from those of $\Phi^A_{\eta, v_k}(z)$, if $j \not= k$.
Then (ii) is proved.
By (\ref{eqn:M2}), 
$\cM^A_{\eta, v_k}(0, x)=\lim_{t \downarrow 0}
\widetilde{\rE}[\Phi^A_{\eta, v_k}(x+i \widetilde{W}(t))]
=\Phi^A_{\eta, v_k}(x), 1 \leq k \leq N$.
Since $\Phi^A_{\eta, v_k}(v_j)=\delta_{jk}, 1 \leq j, k \leq N$
by definition (\ref{eqn:Phi1}),
(iii) is also satisfied. \qed
\vskip 0.3cm

Let $\widetilde{W}_j(\cdot), 1 \leq j \leq N$
be independent $N$ copies of $\widetilde{W}(\cdot)$.
For $\widetilde{\bW}(t)=(\widetilde{W}_1(t), \dots, \widetilde{W}_N(t)), t \geq 0$,
the probability space is denoted by
$(\widetilde{\Omega}_W, \widetilde{\cF}_W, \widetilde{\rP})$
with expectation $\widetilde{\rE}$.
We put
$$
Z^r_j(t)=W^r_j(t)+i \widetilde{W}_j(t), \quad 1 \leq j \leq N,
\quad t \geq 0,
$$
which are independent complex Brownian motions on 
$\C(r) \equiv [0, 2 \pi r) \times i \R$.
The probability space for $\bZ^r(t)=(Z^r_1(t), \dots, Z^r_N(t))$,
$t \geq 0$ is given by the direct product of the two
spaces, $(\Omega_{W^r}, \cF_{W^r}, \rP_{\u}^r)$ for $\bW^r(\cdot)$
and $(\widetilde{\Omega}_W, \widetilde{\cF}_W, \widetilde{\rP})$
for $\widetilde{\bW}(\cdot)$, which is denoted by
$(\Omega^r, \cF^r, \bP_{\u}^r)$ with expectation $\bE_{\u}^r$.

\begin{prop}
\label{thm:expectation1}
Let $\overline{Z}^r_{\delta}(t)=\delta+\sum_{j=1}^N Z^r_j(t)$
and $\overline{W}^r_{\delta}=\delta+\sum_{j=1}^N W^r(t), t \geq 0$.
Then the following equality holds,
\begin{eqnarray}
&& \widetilde{\rE} \left[
\vartheta_1 \left(\frac{\overline{Z}^r_{\delta}(t)}{2 \pi r}; 
\frac{iN t_{*}}{2 \pi r^2} \right)
\prod_{1 \leq j < k \leq N}
\vartheta_1 \left( \frac{Z^r_j(t)-Z^r_k(t)}{2 \pi r}; 
\frac{iN t_{*}}{2 \pi r^2} \right)
\right] 
\nonumber\\
&& = \left[ e^{N t/24 r^2}
\prod_{n=1}^{\infty}
\left( \frac{1-e^{-nN(t_{*}-t)/r^2}}{1-e^{-nN t_{*}/r^2}} \right)
\right]^{-(N-1)(N-2)/2}
\nonumber\\
&&  \times 
\vartheta_1 \left( \frac{\overline{W}^{r}_{\delta}(t)}{2 \pi r}; 
\frac{i N (t_{*}-t)}{2 \pi r^2} \right)
\prod_{1 \leq j < k \leq N}
\vartheta_1 \left(
\frac{W^r_j(t)-W^r_k(t)}{2 \pi r};
\frac{i N (t_{*}-t)}{2 \pi r^2} \right).
\nonumber
\end{eqnarray}
\end{prop}
\noindent{\it Proof.} 
By (\ref{eqn:equalityB2}) in Lemma \ref{thm:equalityB2}
\begin{eqnarray}
&& 
\widetilde{\rE} \left[
\vartheta_1 \left( 
\frac{\overline{Z}^r_{\delta}(t)}{2 \pi r}; \tau \right)
\prod_{1 \leq j < k \leq N}
\vartheta_1 \left( \frac{Z^r_j(t)-Z^r_k(t)}{2 \pi r}; \tau \right)
\right] 
\nonumber\\
&=& 
C_N^{A}(\tau)
\det_{1 \leq j, k \leq N}
\left[ 
\widetilde{\rE} \left[
e^{i(k-1) Z^r_j(t)/r}
\vartheta_1 \left(
\frac{N-1}{2}+(k-1)\tau+ \frac{\delta+N Z^r_j(t)}{2 \pi r}; N \tau \right) \right] \right],
\label{eqn:AA1}
\end{eqnarray}
where the multilinearity of determinant
and independence of $Z^r_j(t)$'s have been used.
Using the Laurent expansion (\ref{eqn:Jacobi3}), we have
\begin{eqnarray}
&&
\widetilde{\rE} \left[
e^{i (k-1) Z^r_j(t)/r}
\vartheta_1 \left(
\frac{N-1}{2}+(k-1)\tau+ \frac{\delta+N Z^r_j(t)}{2 \pi r}; N \tau \right) \right]
\nonumber\\
&&
=
e^{i(k-1) W^r_j(t)/r} 
i \sum_{n \in \Z} (-1)^n 
e^{[(n-(1/2))^2N \tau +(2n-1)\{(N-1)/2+(k-1)\tau+(\delta+N W^r_j(t))/2 \pi r\}] \pi i}
\nonumber\\
&& \qquad 
\times \widetilde{\rE}
\Big[ e^{-\{2(k-1)+(2n-1)N\}\widetilde{W}_j(t)/2r} \Big].
\label{eqn:AA2}
\end{eqnarray}
Here
\begin{eqnarray}
&& \widetilde{\rE} \Big[ e^{-\{2(k-1)+(2n-1)N\}\widetilde{W}_j(t)/2r} \Big]
= \int_{\R} d \widetilde{w} \, \frac{e^{-\widetilde{w}^2/2t}}{\sqrt{2 \pi t}}
e^{-\{2(k-1)+(2n-1)N\} \widetilde{w}/2r}
\nonumber\\
&& \quad = e^{t[2(k-1)+(2n-1)N]^2/8r^2}
\nonumber\\
&& \quad = e^{(k-1)^2 t/2r^2}
\exp \left[
\left(n-\frac{1}{2}\right)^2 N \left(- \frac{i N t}{2 \pi r^2} \right) \pi i
+(2n-1)(k-1) \left( - \frac{i N t}{2 \pi r^2} \right) \pi i \right].
\nonumber
\end{eqnarray}
Then (\ref{eqn:AA2}) is equal to
\begin{eqnarray}
&& e^{(k-1)^2 t/2r^2}
e^{i(k-1)W^r_j(t)/r}
i \sum_{n \in \Z} (-1)^n
\exp \left[ \left(n-\frac{1}{2} \right)^2 N 
\left(\tau- \frac{i N t}{2 \pi r^2} \right) \pi i
\right.
\nonumber\\
&& \qquad \qquad  \left.
+(2n-1) \left\{ \frac{N-1}{2}+(k-1) \left(\tau- \frac{i N t}{2 \pi r^2} \right)
+ \frac{\delta + N W^r_j(t)}{2 \pi r} \right\} \pi i 
\right]
\nonumber\\
&=& e^{(k-1)^2 t/2r^2}
e^{i(k-1)W^r_j(t)/r}
\nonumber\\
&& \times
\vartheta_1 \left(
\frac{N-1}{2}+(k-1) \left(\tau- \frac{i N t}{2 \pi r^2} \right) + 
\frac{\delta+N W^r_j(t)}{2 \pi r}; N \left(\tau-\frac{i N t}{2 \pi r^2} \right) \right).
\nonumber
\end{eqnarray}
Put this into (\ref{eqn:AA1}), we have
\begin{eqnarray}
&& 
\widetilde{\rE} \left[
\vartheta_1 \left( \frac{\overline{Z}^r_{\delta}(t)}{2 \pi r}; \tau \right)
\prod_{1 \leq j < k \leq N}
\vartheta_1 \left( \frac{Z^r_j(t)-Z^r_k(t)}{2 \pi r}; \tau \right) \right]
\nonumber\\
&&
= C_N^{A}(\tau) e^{t \sum_{k=1}^N (k-1)^2/2r^2} 
\nonumber\\
&& \, \times
\det_{1 \leq j, k \leq N}
\left[ e^{i(k-1) W^r_j(t)/r}
\vartheta_1 \left( \frac{N-1}{2}+(k-1) \left(\tau- \frac{i N t}{2\pi r^2} \right) 
+\frac{\delta+ N W^r_j(t)}{2 \pi r}; 
N\left(\tau- \frac{i N t}{2 \pi r^2} \right) \right) \right]
\nonumber\\
&&
= \frac{C_N^{A}(\tau)}{C_N^{A}(\tau-i N t/2 \pi r^2)}
e^{(N-1)N(2N-1) t/12r^2}
C_N^{A} \left(\tau- \frac{i N t}{2 \pi r^2} \right)
\nonumber\\
&& \, \times 
\det_{1 \leq j, k \leq N}
\left[ e^{i(k-1) W^r_j(t)/r}
\vartheta_1 \left( \frac{N-1}{2}+(k-1) \left(\tau-\frac{i N t}{2\pi r^2} \right) 
+ \frac{\delta+ N W^r_j(t)}{2 \pi r}; N \left(\tau- \frac{i N t}{2\pi r^2} \right) \right) \right]
\nonumber\\
&&
= \frac{C_N^{A}(\tau)}{C_N^{A}(\tau-i N t/2 \pi r^2)}
e^{(N-1)N(2N-1) t/12r^2}
\nonumber\\
&& \, \times
\vartheta_1 \left( \frac{\overline{W}^r_{\delta}(t)}{2 \pi r}; 
\tau- \frac{i N t}{2 \pi r^2} \right)
\prod_{1 \leq j < k \leq N}
\vartheta_1 \left( \frac{W^r_j(t)-W^r_k(t)}{2 \pi r}; 
\tau- \frac{i N t}{2 \pi r^2} \right),
\nonumber
\end{eqnarray}
where (\ref{eqn:equalityB2}) of Lemma \ref{thm:equalityB2} was used again.
By (\ref{eqn:equalityB2b}), 
\begin{eqnarray}
\frac{C_N^{A}(\tau)}{C_N^{A}(\tau-i N t/2 \pi r^2)}
&=& e^{-(N-1)N(3N-2) t/16r^2}
\left(\frac{q_0(\tau)}{q_0(\tau-iNt/2 \pi r^2)}
\right)^{(N-1)(N-2)/2}
\nonumber\\
&=& e^{-(N-1)N(3N-2) t/16r^2}
\prod_{n=1}^{\infty}
\left( \frac{1-e^{2 n \pi i \tau}}{1-e^{2 n \pi i \tau +nNt/r^2}} \right)^{(N-1)(N-2)/2}.
\nonumber
\end{eqnarray}
If we set $\tau=iNt_{*}/2 \pi r^2$, 
the equality is obtained. \qed
\vskip 0.3cm

For $\bW^r(t), t \geq 0$, define
\begin{equation}
\cD^A_{\xi}(t, \bW^r(t))
=\det_{1 \leq j, k \leq N}
[\cM^A_{\xi, u_k}(t, W^r_j(t))],
\quad t \in [0, t_{*}),
\label{eqn:D1}
\end{equation}
which we call the {\it determinantal martingale} \cite{Kat13a}.
By Lemma \ref{thm:martingale}, it is a continuous-time martingale.
Then the following equality is established.

\begin{lem}
\label{thm:det_h}
Assume that $\xi=\sum_{j=1}^N u_j \in \mM_0([0, 2 \pi r))$
and $\overline{u}_{\delta} \in (0, 2 \pi r)$. Then
$$
\cD^A_{\xi}(t, \bW^r(t))
= \frac{h_N^{A}(t_{*}-t, \bW^r(t))}{h_N^{A}(t_{*}, \u)},
\quad t \in [0, t_{*}).
$$
\end{lem}
\noindent{\it Proof.} 
By multilinearity of determinant and
independence of $\widetilde{W}_j(\cdot), 1 \leq j \leq N$,
(\ref{eqn:D1}) with (\ref{eqn:Phi1}) and (\ref{eqn:M2}) gives
\begin{eqnarray}
\cD^A_{\xi}(t, \bW^r(t))
&=& \widetilde{\rE} \Bigg[
\det_{1 \leq j, k \leq N} \Bigg[
\frac{\vartheta_1((\overline{u}_{\delta}+Z^r_j(t)-u_k)/2 \pi r; iNt_{*}/2 \pi r^2)}
{\vartheta_1(\overline{u}_{\delta}/2 \pi r; iNt_{*}/2 \pi r^2)}
\nonumber\\
&& \qquad \qquad \left. \left.
\times \prod_{\substack{1 \leq \ell \leq N, \cr \ell \not=k}}
\frac{\vartheta_1((Z^r_j(t)-u_{\ell})/2 \pi r; iNt_{*}/2 \pi r^2)}{\vartheta_1((u_k-u_{\ell})/2 \pi r; iNt_{*}/2 \pi r^2)}
\right]
\right].
\nonumber
\end{eqnarray}
By Lemma \ref{thm:equalityB1}, it is equal to
$$
\widetilde{\rE} \left[
\frac{\vartheta_1(\overline{Z}^r_{\delta}(t)/2 \pi r; iNt_{*}/2 \pi r^2)}
{\vartheta_1(\overline{u}_{\delta}/2 \pi r; iNt_{*}/2 \pi r^2)}
\prod_{1 \leq j < k \leq N}
\frac{\vartheta_1((Z^r_j(t)-Z^r_k(t))/2 \pi r; iNt_{*}/2 \pi r^2)}
{\vartheta_1((u_j-u_k)/2 \pi r; iNt_{*}/2 \pi r^2)}
\right].
$$
Then we apply Proposition \ref{thm:expectation1}.
By definition (\ref{eqn:h1}) of $h_N^{A}$, the equality is obtained. \qed

\SSC{Main Results \label{sec:results}}
\subsection{Determinantal martingale representation \label{sec:DMR}}
By Lemmas \ref{thm:det_h} we obtain the
following representation.
We call it the {\it determinantal martingale representation} (DMR)
for the process $(\Xi^A(t), t \in [0, t_{*}), \P^A_{\eta})$.
\begin{thm}
\label{thm:DMR}
Suppose that $N \in \N$,
$\eta=\sum_{j=1}^{N} \delta_{v_j}$
with (\ref{eqn:ed1}) and (\ref{eqn:delta_A}).
Let $T \in [0, t_{*})$.
For any $\cF_{\Xi^A}(T)$-measurable observable $F$, 
\begin{eqnarray}
\E^A_{\eta} \left[ F \left(\Xi^A(\cdot) \right) \right]
&=& \rE_{\v}^r \left[F \left( \sum_{j=1}^{N} \delta_{W^r_j(\cdot)} \right)
\cD^A_{\eta}(T, \bW^r(T)) \right]
\nonumber\\
&=& \bE_{\v}^r \left[F \left( \sum_{j=1}^{N} \delta_{\Re Z^r_j(\cdot)} \right)
\det_{1 \leq j, k \leq N}
[\Phi^A_{\eta, v_k}(Z^r_j(T))] \right]. 
\label{eqn:DMR1}
\end{eqnarray}
\end{thm}
\vskip 0.3cm
Note that the second representation of (\ref{eqn:DMR1}) 
is an elliptic extension of the 
{\it complex Brownian motion representation}
reported in \cite{KT13} for the Dyson model
({\it i.e.} the noncolliding BM).
\vskip 0.3cm
\noindent{\it Proof.} 
It is sufficient to consider the case that 
$F$ is given as (\ref{eqn:measure1}).
Moreover, by Markov property, it is enough to prove the case
$M=1$; $0 \leq t_1 \leq T < \infty$.
Here we prove the equalities
\begin{eqnarray}
\E^A_{\eta}\left[ g_1(\X^A(t_1)) \right]
&=& \rE_{\v}^r \left[ g_1(\bW^r(t_1))
\cD^A_{\eta}(t_1, \bW^r(t_1)) \right]
\nonumber\\
&=& \bE_{\v}^r \left[g_1(\bW^r(t_1))
\det_{1 \leq j, k \leq N}
[\Phi^A_{\eta, v_k}(Z^r_j(t_1))] \right],
\label{eqn:DMRa1}
\end{eqnarray}
where $g_1$ is a symmetric function
having periodicity (\ref{eqn:g_periodic}). 
By Proposition \ref{thm:h_trans2}, 
\begin{equation}
\E^A_{\eta}\left[ g_1(\X^A(t_1)) \right]
=\check{\rE}_{\v} \left[g_1(\check{\bW}(t_1))
\1(T_{\check{\bbW}} \wedge T_{\overline{W}_{\delta}} > t_1)
\frac{h_N^{A}(t_{*}-t_1, \check{\bW}(t_1))}{h_N^{A}(t_{*}, \v)} \right].
\label{eqn:h_trans2b}
\end{equation}
The definition of $h_N^A$ given by (\ref{eqn:h1})
and the initial condition $\eta$ give 
$$
\mbox{(RHS)}
= \check{\rE}_{\v} \left[g_1(\check{\bW}(t_1))
\1(T_{\check{\bbW}} \wedge T_{\overline{W}_{\delta}} > t_1)
\frac{|h_N^{A}(t_{*}-t_1, \check{\bW}(t_1))|}{h_N^{A}(t_{*}, \v)} \right].
$$
By the determinantal formula (\ref{eqn:qNA_1})
of $q_N^A$ given in Proposition \ref{thm:qNA},
the above is written as
\begin{equation}
\rE_{\v}^r \left[
\sum_{\sigma \in \cS_N} {\rm sgn}(\sigma)
g_1(\bW^r(t_1))
\1(\sigma(\bW^r(t_1)) \in \cA_{[0, 2 \pi r)^N})
\frac{|h_N^{A}(t_{*}-t_1, \bW^r(t_1))|}{h_N^{A}(t_{*}, \v)} \right],
\label{eqn:Esgn}
\end{equation}
where $\bW^r(t_1)=(W^r_1(t_1), \dots, W^r_N(t_1))$ and 
the transition density of 
each $W^r_j$ is given by (\ref{eqn:prZ1}), $1 \leq j \leq N$.
Here we used the notation
$\sigma(\x)=(x_{\sigma(1)}, \dots, x_{\sigma(N)})$ for $\sigma \in \cS_N$.
Since
\begin{eqnarray}
&&
{\rm sgn}(\sigma) \1(\sigma(\bW^r(t_1)) \in \cA_{[0, 2 \pi r)^N})
|h^A_N(t_*-t_1, \bW^r(t_1))|
\nonumber\\
&& \quad 
= \1(\sigma(\bW^r(t_1)) \in \cA_{[0, 2 \pi r)^N})
h^A_N(t_*-t_1, \bW^r(t_1)), \quad \sigma \in \cS_N, 
\nonumber
\end{eqnarray}
(\ref{eqn:Esgn}) is equal to
\begin{eqnarray}
&&
\rE_{\v}^r \left[
\sum_{\sigma \in \cS_N} \1(\sigma(\bW^r(t_1)) \in \cA_{[0, 2 \pi r)^N})
g_1(\bW^r(t_1)) \frac{h^A_N(t_*-t_1, \bW^r(t_1))}{h^A_N(t_*, \u)} \right]
\nonumber\\
&& \quad 
= \rE_{\v}^r \left[ g_1(\bW^r(t_1))
\frac{h^A_N(t_*-t_1, \bW^r(t_1))}{h^A_N(t_*, \u)} \right].
\nonumber
\end{eqnarray}
Then by Lemma \ref{thm:det_h}, 
we obtain the first line of (\ref{eqn:DMRa1}). 
By definitions of $\bE_{\v}^r$ and $\cD^A_{\eta}$ given by (\ref{eqn:D1})
with (\ref{eqn:M2}), the second line 
of (\ref{eqn:DMRa1}) is also obtained. 
\qed
\vskip 0.3cm
\noindent{\bf Remark 2.} \,
The function $h^A_N(t_*-t, \x), t \in [0, t_*)$
is not a harmonic function of $\x$, but
Lemma \ref{thm:martingale} proves that
$\cD^A_{\eta}(t, \bW^r(t))$ given by (\ref{eqn:D1})
is a continuous-time martingale, where $\bW^r(t)$ is a Markov process
defined by using Brownian motion in Section \ref{sec:bWr}.
Then It\^o's formula implies 
$$
\left( \frac{\partial}{\partial t} + \frac{1}{2} \Delta \right)
\cD^A_{\eta}(t, \x)=0,
$$
where $\Delta=\sum_{j=1}^N \partial^2/\partial x_j^2$.
In this sense, DMR is a time-dependent extension
of $h$-transform \cite{Kat13a}.
\vskip 0.3cm

\subsection{Determinantal process \label{sec:det_process}}
For any integer $M \in \N$,
a sequence of times
$\t=(t_1,\dots,t_M)$ with 
$0 \leq t_1 < \cdots < t_M < t_{*}$,
and a sequence of functions
$\f=(f_{t_1},\dots,f_{t_M}) \in \rC([0, 2 \pi r))^M$,
the {\it moment generating function} of multitime distribution
of $(\Xi^A(t), t \in [0,t_{*}), \P^A_{\xi})$ is defined by
\begin{equation}
\Psi^A_{\xi, \t}[\f]
= \E^A_{\xi} \left[ \exp \left\{ \sum_{m=1}^{M} 
\int_{0}^{2 \pi r} f_{t_m}(x) \Xi(t_m, dx) \right\} \right].
\label{eqn:GF1}
\end{equation}
It is expanded with respect to `test functions' 
$\chi_{t_m}(\cdot)=e^{f_{t_m}(\cdot)}-1, 
1 \leq m \leq M$
as
$$
\Psi^A_{\xi, \t}[\f]
=\sum_
{\substack
{0 \leq N_m \leq N, \\ 1 \leq m \leq M} }
\int_{\prod_{m=1}^{M} \cA_{[0, 2 \pi r)^{N_m}}}
\prod_{m=1}^{M} \left\{ d \x_{N_m}^{(m)}
\prod_{j=1}^{N_{m}} 
\chi_{t_m} \Big(x_{j}^{(m)} \Big) \right\}
\rho_{\xi} 
\Big( t_{1}, \x^{(1)}_{N_1}; \dots ; t_{M}, \x^{(M)}_{N_M} \Big),
$$
and it defines the {\it spatio-temporal correlation functions}
$\rho_{\xi}(\cdot)$ for the process $(\Xi^A(t), t \in [0,t_{*}), \P^A_{\xi})$.

Given an integral kernel
$
\mbK(s,x;t,y), 
(s,x), (t,y) \in [0, t_{*}) \times [0, 2 \pi r),
$
the {\it Fredholm determinant} is defined as
\begin{eqnarray}
&& \mathop{{\rm Det}}_
{\substack{
(s,t)\in \{t_1, \dots, t_M\}^2, \\
(x,y)\in [0, 2 \pi r)^2}
}
 \Big[\delta_{st} \delta_x(y)
+ \mbK(s,x;t,y) \chi_{t}(y) \Big]
\nonumber\\
&& 
=\sum_
{\substack
{0 \leq N_m \leq N, \\ 1 \leq m \leq M} }
\sum_
{\substack
{\x^{(m)}_{N_m} \in \cA_{[0, 2 \pi r)^{N_m}}, 
\\ 1 \leq m \leq M} }
\prod_{m=1}^{M}
\prod_{j=1}^{N_{m}} 
\chi_{t_m} \Big(x_{j}^{(m)} \Big)
\det_{\substack
{1 \leq j \leq N_{m}, 1 \leq k \leq N_{n}, \\
1 \leq m, n \leq M}
}
\Bigg[
\mbK(t_m, x_{j}^{(m)}; t_n, x_{k}^{(n)} )
\Bigg].
\label{eqn:F_det}
\end{eqnarray}
We put the following definition \cite{BR05,KT10}.
\begin{df}
\label{thm:determinantal}
For a given initial configuration $\xi$, 
if any moment generating function (\ref{eqn:GF1}) 
is expressed by a Fredholm determinant, we say
the process $(\Xi^A(t), t \in [0, t_{*}), \P^A_{\xi})$ is determinantal. 
In this case, 
all spatio-temporal correlation functions
are given by determinants as 
\begin{equation}
\rho_{\xi} \Big(t_1,\x^{(1)}_{N_1}; \dots;t_M,\x^{(M)}_{N_M} \Big) 
=\det_{\substack
{1 \leq j \leq N_{m}, 1 \leq k \leq N_{n}, \\
1 \leq m, n \leq M}
}
\Bigg[
\mbK_{\xi}(t_m, x_{j}^{(m)}; t_n, x_{k}^{(n)} )
\Bigg],
\label{eqn:rho1}
\end{equation}
$0 \leq t_1 < \cdots < t_M < t_{*}$,
$1 \leq m \leq M$, 
$1 \leq N_m \leq N$,
$\x^{(m)}_{N_m} \in [0, 2 \pi r)^{N_m}, 1 \leq m \leq M \in \N$.
Here the integral kernel $\mbK_{\xi}: ([0, t_{*}) \times [0, 2 \pi r))^2 \mapsto \R$ is called
the (spatio-temporal) correlation kernel.
\end{df}

By Theorem 1.3 in \cite{Kat13a},
DMR given by Theorem \ref{thm:DMR}
leads to the following result.

\begin{cor}
\label{thm:kernel}
For $\eta=\sum_{j=1}^N \delta_{v_j}$ with
(\ref{eqn:ed1}) and (\ref{eqn:delta_A}), 
the process $(\Xi^A(t), t \in [0, t_{*}), \P^A_{\eta})$ is
determinantal with the correlation kernel
\begin{eqnarray}
\mbK^A_{\eta}(s, x; t, y)
&=& \mbK^A_{\eta}(s,x;t,y; N, r, t_{*})
\nonumber\\
&=& \int_0^{2 \pi r} \eta(du) \,
p^r_{A_{N-1}}(s,x|u) \cM^A_{\eta, u}(t,y)-\1(s>t) p^r_{A_{N-1}}(s-t,x|y),
\label{eqn:K1}
\end{eqnarray}
$(s,x), (t,y) \in [0, t_{*}) \times [0, 2 \pi r)$.\\
\end{cor}

\subsection{Explicit expression of $\mbK^A_{\eta}$ 
and infinite-particle limit \label{sec:special_initial}}

For $\eta=\sum_{j=1}^N \delta_{v_j}$ with
(\ref{eqn:ed1}) and (\ref{eqn:delta_A}), 
the entire functions (\ref{eqn:Phi1}) become
\begin{eqnarray}
\Phi^A_{\eta, v_k}(z; N, r, t_{*})
&=& \frac{\vartheta_1(z/2\pi r-(k-1)/N +1/2; \tau)}
{\vartheta_1(1/2; \tau)}
\prod_{\substack{1 \leq \ell \leq N, \cr\ell \not=k}}
\frac{\vartheta_1(z/2\pi r-(\ell-1)/N; \tau)}
{\vartheta_1((k-\ell)/N; \tau)}
\nonumber\\
&=& \frac{\vartheta_1(z/2\pi r-(k-1)/N+1/2; \tau)}
{\vartheta_1(1/2; \tau)}
\prod_{n=1}^{N-1}
\frac{\vartheta_1(z/2 \pi r-(k-1)/N+n/N; \tau)}
{\vartheta_1(n/N; \tau)},
\nonumber
\end{eqnarray}
$1 \leq k \leq N$, 
with $\tau=\tau(t_*)=i N t_*/2 \pi r^2$, 
where we have used (\ref{eqn:Jacobi5a}).
Using the formulas (\ref{eqn:formula1}) and (\ref{eqn:formula2}),
it is written as
\begin{eqnarray}
\Phi^A_{\eta, v_k}(z; N, r, t_{*})
&=& \frac{\pi }{N \vartheta_1'(0; N \tau)} \vartheta_1(N\{z/2\pi r- (k-1)/N \}; N \tau)
\nonumber\\
&& \times
\frac{\vartheta_1(z/2 \pi r - (k-1)/N+1/2; \tau)
\vartheta_1'(0; \tau)}
{\pi \vartheta_1(z/2 \pi r-(k-1)/N; \tau) \vartheta_1(1/2; \tau) }, 
\label{eqn:Phi_elliptic}
\end{eqnarray}
$1 \leq k \leq N$.
If we apply the formula (\ref{eqn:formula3}) and the Laurent expansion
(\ref{eqn:Jacobi3}) of $\vartheta_1$, we have
\begin{eqnarray}
&& \Phi^A_{\eta, v_k}(z; N, r, t_{*})
= \frac{2 \pi e^{N \pi i \tau/4}}
{N \vartheta_1'(0; N \tau)}
\nonumber\\
&& \quad \times \left[
\cos(z/2r-(k-1)\pi/N) \sum_{n=1}^{\infty} (-1)^{n-1} e^{N \pi i \tau n(n-1)}
\frac{\sin \Big[(2n-1)N\{z/2r-(k-1)\pi/N\} \Big]}{\sin(z/2r-(k-1)\pi/N)}
\right.
\nonumber\\
&& \qquad \quad -4 \sum_{n=1}^{\infty} 
(-1)^{n-1} e^{N \pi i \tau n(n-1)}  \sum_{\ell=1}^{\infty} 
\frac{e^{2\pi i \tau \ell}}{1+e^{2 \pi i \tau \ell}}
\nonumber\\
&& \qquad \qquad \qquad \qquad 
\times \sin \Big[(2n-1)N\{z/2r-(k-1) \pi/N\} \Big] 
\sin \Big[2 \ell \{z/2r-(k-1) \pi /N\} \Big] \Bigg],
\nonumber
\end{eqnarray}
since $\cot(\pi/2)=0$ and
$\sin(\theta+m \pi)=(-1)^m \sin \theta, m \in \Z$.
For $M \in \N$, let
\begin{equation}
\sigma_M(m)= \left\{
\begin{array}{ll}
m, \quad & \mbox{if $M$ is odd}, \cr
m-1/2, \quad & \mbox{if $M$ is even}.
\end{array} \right.
\label{eqn:spin}
\end{equation}
It is easy to confirm the equality
$$
\frac{\sin(M x)}{\sin x}
=\sum_{\substack{m \in \Z, \cr |\sigma_M(m)| \leq (M-1)/2}}
e^{2 i \sigma_M(m)x}.
$$
Then we see
\begin{eqnarray}
&& \Phi^A_{\eta, v_k}(z; N, r, t_*)=
\frac{2 \pi e^{N \pi i \tau/4}}{N \vartheta_1'(0; N \tau)}
\nonumber\\
&& \times \Bigg[ 
\sum_{n=1}^{\infty}(-1)^{n-1} e^{N \pi i \tau n(n-1)} \Bigg\{
\cos \Big[(2n-1)N\{z/2r-(k-1)\pi/N\} \Big]
\nonumber\\
&& \qquad \qquad \qquad \qquad \qquad \qquad \qquad
+ \sum_{\substack{m \in \Z, \cr
|\sigma_{N-1}(m)| \leq \{(2n-1)N-2\}/2}}
e^{2 i \sigma_{N-1}(m) \{z/2r-(k-1) \pi/N\}} \Bigg\}
\nonumber\\
&& \quad +2 \sum_{n \in \Z}(-1)^{n-1} e^{N \pi i \tau n(n-1)}
\sum_{\ell=1}^{\infty} 
\frac{e^{2 \pi i \tau \ell}}{1+e^{2 \pi i \ell}}
\cos\Big[\{(2n-1)N+2 \ell\}\{z/2r-(k-1)\pi/N\} \Big] \Bigg].
\nonumber
\end{eqnarray}
Here we have used the fact that, for $M \in \N$, 
$2 \sigma_M(m)+1=2 \sigma_{M-1}(m+1)$
if $M$ is odd, and
$2 \sigma_M(m)+1=2 \sigma_{M-1}(m)$ if $M$ is even,
and that $\sigma_{(2n-1)N-1}(m)=\sigma_{N-1}(m)$
for $n, N \in \N$.

The functions 
\begin{eqnarray}
\cM^A_{\eta, v_k}(t, x; N, r, t_*)
&=& \widetilde{\rE}[ \Phi^A_{\eta, v_k}(x+i \widetilde{W}(t); N, r, t_*)]
\nonumber\\
&=& \int_{\R} d \widetilde{w} \,
\frac{e^{-\widetilde{w}^2/2t}}{\sqrt{2 \pi t}}
\Phi^A_{\eta, v_k}(x+i \widetilde{w}; N, r, t_*),
\quad 1 \leq k \leq N, 
\nonumber
\end{eqnarray}
which give continuous-time martingales if we put
$x=W^r_j(\cdot), 1 \leq j \leq N$ (Lemma \ref{thm:martingale} (i)), 
are calculated by performing Gaussian integrals
for each term. 
By setting $\tau=i N t_*/2 \pi r^2$,
the result is expressed as
\begin{eqnarray}
&& \cM^A_{\eta, v_k}(t, x; N, r, t_*)
= \frac{2 \pi}
{N \vartheta_1'(0; i N^2 t_*/2 \pi r^2)}
\nonumber\\
&& \, \times \Bigg[
\sum_{n=1}^{\infty}(-1)^{n-1} e^{-(n-(1/2))^2 N^2(t_*-t)/2r^2} 
\cos \Big[(2n-1)N\{x/2r-(k-1)\pi/N\} \Big]
\nonumber\\
&& \quad
+ \sum_{n=1}^{\infty}(-1)^{n-1} e^{-(n-(1/2))^2 N^2 t_*/2r^2} 
\nonumber\\
&& \qquad \qquad \times
\sum_{\substack{m \in \Z, \cr 
|\sigma_{N-1}(m)| 
\leq \{(2n-1)N-2\}/2}}
e^{\sigma_{N-1}(m)^2 t/2r^2+2i \sigma_{N-1}(m)\{x/2r-(k-1) \pi/N\}}
\nonumber\\
&& \quad
+ 2 \sum_{n \in \Z} (-1)^{n-1} e^{-(n-(1/2))^2 N^2(t_*-t)/2r^2}
\sum_{\ell=1}^{\infty} 
\frac{e^{-\ell N t_*/r^2}}{1+e^{-\ell N t_*/r^2}}
e^{\ell \{\ell+(2n-1)N\}t/2r^2}
\nonumber\\
&& \qquad \qquad \qquad \qquad
\times \cos \Big[\{(2n-1)N+2 \ell\} \{x/2r-(k-1)\pi/N\} \Big] \Bigg].
\label{eqn:M_eta2A}
\end{eqnarray}

Then by Corollary \ref{thm:kernel}, the correlation kernel
is determined as
$$
\mbK^A_{\eta}(s, x; t, y; N, r, t_*)
= \cG^A_{\eta}(s, x; t, y; N, r, t_*)
-\1(s>t) p^r_{A_{N-1}}(s-t, x|y), 
$$
$(s,x), (t,y) \in [0, \infty) \times [0, 2 \pi r)$,
where
$$
\cG^A_{\eta}(s,x;t,y;N,r, t_*)
=\sum_{k=1}^N p^r_{A_{N-1}}(s,x| v_k)
\cM^A_{\eta, v_k}(t, y; N, r, t_*)
$$
with (\ref{eqn:prZ3}).
We note that, by using (\ref{eqn:spin}), 
(\ref{eqn:prZ3}) is written as
\begin{equation}
p^r_{A_{N-1}}(t,y|x) 
= \frac{1}{2 \pi r}
\sum_{\ell \in \Z} e^{-\sigma_{N-1}(\ell)^2 t/2r^2+i \sigma_{N-1}(\ell) (y-x)/r},
\quad t \geq 0, \quad x, y \in [0, 2 \pi r).
\label{eqn:pR2}
\end{equation}
We find that, by using the identity
$\sum_{k=1}^N e^{-2i (k-1) \alpha \pi/N}=N \sum_{k \in \Z} \1(\alpha=k N)$, 
the above is expressed as follows,
\begin{eqnarray}
&& \cG^A_{\eta}(s, x; t, y; N, r, t_*)
=\frac{1}{\vartheta_1'(0; i N^2 t_*/2 \pi r^2) r}
\nonumber\\
&& \quad
\times \Bigg[
\frac{1}{2} \sum_{k \in \Z}
e^{-k^2 N^2 s/2r^2+i k N x/r}
\vartheta_2 \left( \frac{N(y-x)}{2 \pi r}-\frac{i k N^2 s}{2 \pi r^2} ;
\frac{iN^2\{t_*-(t-s)\}}{2 \pi r^2} \right)
\nonumber\\
&& \qquad 
+ \sum_{n=1}^{\infty}(-1)^{n-1} e^{-(n-(1/2))^2 N^2 t_*/2r^2} 
\sum_{\substack{m \in \Z, \cr 
|\sigma_{N-1}(m)| 
\leq \{(2n-1)N-2\}/2}}
e^{\sigma_{N-1}(m)^2 (t-s)/2r^2+i \sigma_{N-1}(m)(y-x)/r}
\nonumber\\
&& \qquad \qquad \qquad \times
\vartheta_3 \left( \frac{Nx}{2 \pi r}-\frac{i N \sigma_{N-1}(m) s}{2 \pi r^2} ;
\frac{iN^2 s}{2 \pi r^2} \right)
\nonumber\\
&& \qquad
+ 2 \sum_{k \in \Z} e^{-k^2 N^2 s/2r^2 + i k N x/r}
\sum_{\ell=1}^{\infty} 
\sum_{n \in \Z} (-1)^{n-1} 
\frac{e^{-\ell N t_*/r^2}}{1+e^{-\ell N t_*/r^2}}
e^{-(n-(1/2))^2 N^2t_*/2r^2}
\nonumber\\
&& \qquad \qquad \qquad 
\times 
e^{\{(n-(1/2))N+\ell\}^2(t-s)/2r^2}
\cos \left[\{(2n-1)N+2 \ell\} \pi
\left( \frac{y-x}{2 \pi r}- \frac{i k N s}{2 \pi r^2} \right) \right],
\label{eqn:GA_1}
\end{eqnarray}
where $\vartheta_2$ and $\vartheta_3$ are defined
by (\ref{eqn:theta_0_2_3}).

We consider the infinite-particle limit with fixed particle-density;
\begin{equation}
N \to \infty, \quad r \to \infty \quad \mbox{with} \quad
\rho = \frac{N}{2 \pi r} = \mbox{const.}
\label{eqn:limitA}
\end{equation}
We can see that the first term in (\ref{eqn:GA_1}) vanishes
in this limit because of the overall factor $1/r$, and obtain
\begin{eqnarray}
&& \bG^A_{\eta}(s, x; t, y; \rho, t_*)
\equiv \lim_{\substack{N \to \infty, r \to \infty, \cr \rho={\rm const.}}}
\cG^A_{\eta}(s, x; t, y; N, r, t_*)
\nonumber\\
&&  \quad
= \frac{2 \pi}{\vartheta_1'(0; 2 \pi i \rho^2 t_*)} \Bigg[
\sum_{n=1}^{\infty} (-1)^{n-1} e^{-2 \pi^2(n-(1/2))^2 \rho^2 t_*}
\nonumber\\
&& \qquad \qquad \qquad \qquad \qquad
\times \int_{|v| \leq (2n-1) \rho} d v \,
e^{\pi^2 v^2(t-s)/2 + i \pi v (y-x)}
\vartheta_3 \Big(\rho x - i \pi v \rho s; 2 \pi i \rho^2 s \Big)
\nonumber\\
&& \qquad \qquad \qquad
+ \sum_{k \in \Z} e^{-2 \pi^2 k^2 \rho^2 s + 2 \pi k \rho x}
\sum_{n \in \Z} (-1)^{n-1}
\int_0^{\infty} dv \,
\frac{e^{-2 \pi^2 v \rho t_*}}{1+e^{-2 \pi^2 v \rho t_*}}
\nonumber\\
&& \qquad \qquad \qquad \quad 
\times e^{\pi^2\{(2n-1) \rho+v\}^2(t-s)/2}
\cos \Big[\{ (2n-1) \rho+v\} \pi (y-x-2 \pi i k \rho s) \Big] \Bigg].
\label{eqn:GNinf1}
\end{eqnarray}
The correlation kernel in this limit (\ref{eqn:limitA})
is given by
\begin{equation}
\bK^A_{\eta}(s,x; t, y; \rho, t_*)
=\bG^A_{\eta}(s, x; t, y; \rho, t_*)
-\1(s>t) p_{\rm BM}(s-t, x|y),
\label{eqn:KNinf1}
\end{equation}
$(s,x), (t,y) \in [0, t_*) \times \R$.
The convergence of correlation kernel (\ref{eqn:GNinf1}) in the limit (\ref{eqn:limitA})
implies well-definedness of elliptic determinantal process
with an infinite number of particles.
Further study will be reported elsewhere
(see Section \ref{sec:reduction_infinite_A} below). 

\SSC{Reduced Processes in Temporally Homogeneous Limit \label{sec:reductions}}
\subsection{Reduction in SDEs \label{sec:reduction_SDE}}
Since we have 
\begin{equation}
\lim_{t_{*} \to \infty} A_N^{2 \pi r}(t_{*}-t,x)
= \frac{1}{2 r} \cot \left( \frac{x}{2r} \right)
\label{eqn:limit1}
\end{equation}
by (\ref{eqn:A3}), 
the limit $t_{*} \to \infty$ of (\ref{eqn:SDEe1}) gives 
a temporally homogeneous system of the SDEs, 
\begin{equation}
d \chX^A_j(t) = dB_j(t)+\frac{1}{2 r}
\sum_{\substack{1 \leq k \leq N, \cr k \not=j}}
\cot \left( \frac{\chX^A_j(t)- \chX^A_k(t)}{2r} \right) dt
+\frac{1}{2 r}
\cot \left( \frac{\overline{X}^A_{\delta}(t)}{2r} \right) dt, 
\quad t \geq 0, 
\label{eqn:circle1}
\end{equation}
$1 \leq j \leq N$, and (\ref{eqn:SDEe2}) becomes
\begin{equation}
d \overline{X}^A_{\delta}(t)
= \sqrt{N} dB(t)+ \frac{N}{2r} \cot \left(
\frac{\overline{X}^A_{\delta}(t)}{2r} \right) dt, 
\quad t \geq 0.
\label{eqn:circle2}
\end{equation}
The system of SDEs without the third term in RHS of (\ref{eqn:circle1})
has been studied as a dynamical extension of the 
{\it circular unitary ensemble (CUE)} of random matrix theory
\cite{HW96,Meh04,For10}.
It is interesting to see that
a one-parameter extension of that system
is discussed as a driving system
for a multiple Schramm-Loewner evolution
by Cardy \cite{Car03}.
Moreover, for
$$
\lim_{r \to \infty} \lim_{t_{*} \to \infty} A_N^{2 \pi r}(t_{*}-t, x) 
= \lim_{r \to \infty} \frac{1}{2r} \cot \left( \frac{x}{2r} \right)
= \frac{1}{x},
$$
and $\delta=\pi r n$ with a fixed $n \in \Z$ determined
by the initial configuration, 
the $r \to \infty$ limit
of the system (\ref{eqn:circle1})  is given by
\begin{equation}
dX^A_j(t)=dB_j(t)+
\sum_{\substack{1 \leq k \leq N, \cr k \not=j}}
\frac{1}{X^A_j(t)-X^A_k(t)} dt, 
\quad 1 \leq j \leq N, \quad t \geq 0,
\label{eqn:DysonSDE1}
\end{equation}
where $X^A_j=\chX^A_j, 1 \leq j \leq N$ in $r \to \infty$.
It is the system of SDEs of the Dyson model
({\it i.e.}, the noncolliding Brownian motion on $\R$).

\subsection{Reduction in correlation kernel of determinantal process \label{sec:reduction_kernel}}
Next we study the reduction as determinantal processes
for $(\Xi^A(t), t \in [0, t_*), \P^A_{\eta})$. 
For the configuration $\eta=\sum_{j=1}^N \delta_{v_j}$
with (\ref{eqn:ed1}) and (\ref{eqn:delta_A}), an explicit expression of
correlation kernel $\mbK^A_{\eta}$ was given in 
Section \ref{sec:special_initial}.
In (\ref{eqn:GA_1}) we can take the temporally homogeneous limit
$t_* \to \infty$ as follows.
We see from (\ref{eqn:theta_dash}) and 
(\ref{eqn:theta_0_2_3}) that
$\vartheta_1'(0; \tau) \sim 2 \pi e^{\pi i \tau/4}$,
$\vartheta_2(v; \tau) \sim 2 e^{\pi i \tau/4} \cos (\pi v)$,
and
$\vartheta_3(v; \tau) \to 1$
as $\Im \tau \to + \infty$.
Then, we have
\begin{eqnarray}
&& \widehat{\cG}^A_{\eta}(s, x; t, y; N, r)
\equiv \lim_{t_* \to \infty} \cG^A_{\eta}(s, x; t, y; N, r, t_*)
\nonumber\\
&& \quad
= \sum_{k \in \Z} e^{-k^2 N^2 s/2r^2 + i k N x/r}
e^{N^2(t-s)/8r^2} \cK^A \left((y-x)- \frac{ikNs}{r}; N, r \right)
\nonumber\\
&& \quad
+ \frac{1}{2 \pi r} \sum_{\substack{m \in \Z, \cr
|\sigma_{N-1}(m)| \leq (N-2)/2}}
e^{\sigma_{N-1}(m)^2 (t-s)/2r^2 + i \sigma_{N-1}(m) (y-x)/r}
\vartheta_3 \left(
\frac{Nx}{2 \pi r} - \frac{i N \sigma_{N-1}(m)s}{2 \pi r^2} ;
\frac{i N^2 s}{2 \pi r} \right),
\nonumber\\
\label{eqn:GA_2}
\end{eqnarray}
where
$$
\cK^A(x;N,r) = \frac{1}{2 \pi r}
\cos \left( \frac{N x}{2 r} \right).
$$
The correlation kernel is given by
\begin{equation}
\widehat{\mbK}^A_{\eta}(s, x; t, y; N, r)
= \widehat{\cG}^A_{\eta}(s, x; t, y; N, r)
-\1(s>t) p^r_{A_{N-1}}(s-t, x|y), 
\label{eqn:K_eta}
\end{equation}
$(s,x), (t,y) \in [0, \infty) \times [0, 2 \pi r)$.
The determinantal process defined by the correlation kernel (\ref{eqn:K_eta})
is denoted by $(\widehat{\Xi}^A(t), t \in [0, \infty), \widehat{\P}^A_{\eta})$.
From the explicit expression (\ref{eqn:GA_2}), we can show that
$(\widehat{\Xi}^A(t), t \in [0, \infty), \widehat{\P}^A_{\eta})$ 
exhibits a typical nonequilibrium phenomenon; 
{\it relaxation to equilibrium}.

Let $(\widehat{\Xi}^A(t), t \in [0, \infty), \widehat{\P}^A_{\rm eq})$
be the equilibrium determinantal process,
whose correlation kernel is homogeneous both in space
$[0, 2 \pi r)$ and time $[0, \infty)$ and given by
\begin{eqnarray}
&& \widehat{\mbK}^A_{\rm eq}(t-s, y-x; N, r)
= \widehat{\cG}^A_{\rm eq}(t-s, y-x; N,r)- \1(s>t) p^r_{A_{N-1}}(s-t, x|y)
\nonumber\\
&& = \left\{ \begin{array}{l}
\displaystyle{ \frac{1}{2\pi r}
\sum_{\substack{\ell \in \Z, \cr |\sigma_{N-1}(\ell)| \leq (N-2)/2}}
e^{\sigma_{N-1}(\ell)^2(t-s)/2r^2+i \sigma_{N-1}(\ell)(y-x)/r}
+e^{N^2(t-s)/8r^2} \cK^A(y-x; N, r)}, \cr
\hskip 8cm \mbox{if $s < t$}, \cr
\displaystyle{
\frac{1}{2 \pi r} \frac{\sin[(N-1)(y-x)/2r]}{\sin[(y-x)/2r]} +\cK^A(y-x; N, r) }, \cr
\hskip 8cm \mbox{if $s=t$}, \cr
\displaystyle{ -\frac{1}{2\pi r}
\sum_{\substack{\ell \in \Z, \cr |\sigma_{N-1}(\ell)| > (N-2)/2}}
e^{\sigma_{N-1}(\ell)^2(t-s)/2r^2+i \sigma_{N-1}(\ell) (y-x)/r} +e^{N^2(t-s)/8r^2} \cK^A(y-x; N, r)}, \cr
\hskip 8cm \mbox{if $s > t$},
\end{array} \right.
\label{eqn:KK1}
\end{eqnarray}
$(s, x), (t,y) \in [0, \infty) \times [0, 2 \pi r)$.
Note that the spatial dependence of (\ref{eqn:KK1})
is on a distance of two points $|y-x|$.

\begin{prop}
\label{thm:relax}
The process $(\widehat{\Xi}^A(t+T), t \in [0, \infty), \widehat{\P}^A_{\eta})$
converges to 
$(\widehat{\Xi}^A(t), t \in [0, \infty), \widehat{\P}^A_{\rm eq})$
as $T \to \infty$ weakly in the sense of finite-dimensional distributions.
\end{prop}
\vskip 0.5cm
Note that if we ignore the terms containing $\cK^A$,
the kernel (\ref{eqn:KK1}) is equal to the correlation kernel
of the equilibrium determinantal process
of the noncolliding BM on a circle $\rS^1(r)$
with $N-1$ particles,
$\mbK_{\rm eq}^{\rm CUE}(t-s, y-x; N-1,r)$ \cite{NF03,Kat13a}.
It is a reversible process with respect to
the eigenvalue-angle distribution
of unitary random matrices with size $N-1$ in the CUE
(if we set $r=1$).
With contribution from the term $\cK^A$,
we have the particle density
\begin{eqnarray}
\rho = \rho(N, r)
&=& \lim_{|y-x| \to 0} \widehat{\mbK}^A_{\rm eq}(0, y-x; N, r)
\nonumber\\
&=& \lim_{|y-x| \to 0} \mbK_{\rm eq}^{\rm CUE}(0, y-x; N-1, r)
+\lim_{|y-x| \to 0} \cK^A(y-x; N, r)
\nonumber\\
&=& \frac{N-1}{2 \pi r}+ \frac{1}{2 \pi r}=\frac{N}{2 \pi r},
\nonumber
\end{eqnarray}
which is uniform on $[0, 2 \pi r)$ as it should be.

\vskip 0.5cm
\noindent{\it Proof of Proposition \ref{thm:relax}.} 
Since $e^{-k^2 N^2 T/2r^2} \to 0, k \not= 0$ and
$\vartheta_3(v; i N^2 T/2 \pi r^2) \to 1$ as $T \to \infty$,
(\ref{eqn:GA_2}) gives
\begin{eqnarray}
&& \lim_{T \to \infty} 
\widehat{\cG}^A_{\eta}(s+T, x; t+T, y; N, r)
\nonumber\\
&& \quad
= e^{N^2(t-s)/8r^2} \cK^A(y-x; N, r)
+ \frac{1}{2 \pi r} \sum_{\substack{\ell \in \Z, \cr |\sigma_{N-1}(\ell)| \leq (N-2)/2}}
e^{\sigma_{N-1}(\ell)^2(t-s)/2r^2+i \sigma_{N-1}(\ell)(y-x)/r}, 
\nonumber
\end{eqnarray}
which is $\widehat{\cG}^A_{\rm eq}(t-s, y-x; N, r)$. 
If we set $t=s$, the second term becomes
$$
 \frac{1}{2 \pi r} \sum_{\substack{\ell \in \Z, \cr |\sigma_{N-1}(\ell)| \leq (N-2)/2}}
e^{i \sigma_{N-1}(\ell)(y-x)/r}
= \frac{1}{2 \pi r} \frac{\sin[(N-1)(y-x)/2r]}{\sin[(y-x)/2r]}.
$$
Combining the above results with the expression (\ref{eqn:pR2}) of $p^r_{A_{N-1}}$,
we obtain (\ref{eqn:KK1}).
The convergence of correlation kernel
$\widehat{\mbK}^A_{\eta} \to \widehat{\mbK}^A_{\rm eq}$ in the long-term limit
guarantees the convergence of 
moment generating function of multitime distribution 
$\Psi^A_{\eta, \t}[\f]$
given by the Fredholm determinant  (\ref{eqn:F_det}) of $\widehat{\mbK}^A_{\eta}$.
It implies the convergence $\widehat{\P}^A_{\eta} \to \widehat{\P}^A_{\rm eq}$
in the long-term limit in the sense of finite-dimensional distributions.
Then the proof is completed. \qed
\vskip 0.5cm

\subsection{Reduction in infinite-particle system \label{sec:reduction_infinite_A}}

The temporally homogeneous limit $t_* \to \infty$
of (\ref{eqn:KNinf1}) with (\ref{eqn:GNinf1}) gives
\begin{eqnarray}
&& \widehat{\bK}^A_{\eta}(s, x; t, y; \rho)
\equiv \lim_{t_* \to \infty} \bK^A_{\eta}(s,x; t, y; \rho, t_*)
\nonumber\\
&& \quad 
= \int_{|v| \leq \rho} d v \,
e^{\pi^2 v^2(t-s)/2+i \pi v(y-x)}
\vartheta_3 \Big(\rho x-i \pi v \rho s; 2 \pi i \rho^2 s \Big)
-\1(s>t) p_{\rm BM}(s-t, x|y),
\label{eqn:Kinf_homo1}
\end{eqnarray}
$(s, x), (t, y) \in [0, \infty) \times \R$.
If we set $\rho=1$, 
this correlation kernel is exactly the same as Eq.(1.5) in \cite{KT10}.
There we gave a set $\mX$ of infinite-particle configurations, 
started at which the Dyson model is well defined as a determinantal process.
The function (\ref{eqn:Kinf_homo1}) was derived as the correlation kernel
of the Dyson model started at $\xi^{\Z} \equiv \sum_{j \in \Z} \delta_{j} \in \mX$.
As shown in \cite{KT10}, for
$(s,x), (t, y) \in [0, \infty) \times \R$,
\begin{eqnarray}
&& \lim_{T \to \infty}
\widehat{\bK}^A_{\eta}(s+T, x; t+T, y; \rho)
\nonumber\\
&&  \qquad 
=\bK_{\rm sin}(t-s, y-x; \rho)
\equiv \left\{ \begin{array}{ll}
\displaystyle{ \int_{0}^{\rho} dv \,
e^{\pi^2 v^2(t-s)/2} \cos[\pi v(y-x)]},
& \mbox{if $s<t$}, \cr
& \cr
\displaystyle{ \frac{\sin[\pi \rho(y-x)]}{\pi(y-x)}},
& \mbox{if $s=t$},\cr
& \cr
\displaystyle{ -\int_{\rho}^{\infty} dv \,
e^{\pi^2 v^2(t-s)/2} \cos[\pi v(y-x)]},
& \mbox{if $s>t$}. \cr
\end{array} \right.
\nonumber
\end{eqnarray}
The limit kernel is 
the {\it extended sine kernel} with density $\rho$ \cite{For10}.
The convergence implies the infinite-dimensional
relaxation phenomenon of the Dyson model
from $\xi^{\Z}$ to equilibrium \cite{KT10}.
This observation is consistent with
Proposition \ref{thm:relax}, since
$\lim_{N \to \infty, r \to \infty, 
\rho={\rm const.}} \cK^A(x; N, r)=0$.

\SSC{Expressions by Gosper's $q$-sine function and hyperbolic limit \label{sec:sinq}}

Gosper defined his $q$-sine function as \cite{Gos01}
$$
\sin_q(\pi z)=q^{(z-1/2)^2}
\frac{(q^{2z}; q^2)_{\infty} (q^{2-2z}; q^2)_{\infty}}{(q; q^2)_{\infty}^2},
\quad 0 < q < 1.
$$
By the product form (\ref{eqn:theta_p1}) of $\vartheta_1(v;\tau)$,
we find the equality
$$
\vartheta_1(v;\tau)
=-i (q;q^2)_{\infty}^2 (q^2;q^2)_{\infty} e^{-2 \pi i v}
\sin_q(-\pi v/\tau)
$$
with (\ref{eqn:ztau1}).
Then the function $A_N^{\alpha}(t_{*}-t,x)$ defined by (\ref{eqn:A2})
is expressed as
\begin{equation}
A_N^{\alpha}(t_{*}-t, x)
= \frac{i \alpha}{2N(t_{*}-t)}
\left. \frac{d}{dz} \log( \sin_q z) \right|_{z=i \alpha x/2N(t_{*}-t)}
-\frac{2 \pi i}{\alpha},
\quad t \in [0, t_{*}), \quad x \in \R, 
\label{eqn:Ab1}
\end{equation}
with
\begin{equation}
q=e^{-2 \pi^2 N(t_{*}-t)/\alpha^2}.
\label{eqn:Ab2}
\end{equation}
The entire functions (\ref{eqn:Phi1}), with which
the determinantal martingales and correlation functions
are expressed, are rewritten as
\begin{eqnarray}
\Phi^A_{\xi, u_k}(z; N, r, t_{*})
&=& e^{-i N(z-u_k)/r}
\frac{\sin_q(\pi i (\overline{u}_{\delta}+z-u_k)r/N t_{*})}
{\sin_q(\pi i \overline{u}_{\delta} r/N t_{*})}
\nonumber\\
&& \times 
\prod_{\substack{1 \leq \ell \leq N, \cr \ell \not= k}}
\frac{\sin_q(\pi i (z-u_{\ell}) r/N t_{*})}{\sin_q(\pi i (u_k-u_{\ell}) r/N t_{*})},
\label{eqn:Phib1}
\end{eqnarray}
$1 \leq k \leq N, z \in \C$,
with $q=e^{-N(t_*-t)/2r^2}$.
We remark that, the $q$-extension of the gamma function
is defined as \cite{AAR99} 
$$
\Gamma_q(z)=(1-q)^{1-z} \frac{(q;q)_{\infty}}{(q^z;q)_{\infty}},
\quad 0 < q < 1,
$$
and the $q$-analogue of Euler's reflection formula
\begin{equation}
\sin_q(\pi z) = q^{1/4} \Gamma_{q^2}(1/2)^2
\frac{(q^2)^{{z \choose 2}}}{\Gamma_{q^2}(z) \Gamma_{q^2}(1-z)}
\label{eqn:gammaq2}
\end{equation}
holds. Then (\ref{eqn:Ab1}) and (\ref{eqn:Phib1}) are
also expressed using the $q$-gamma function
$\Gamma_{q^2}$. 

In the previous section, we studied the temporally homogeneous limit
$t_{*} \to \infty$ with fixed $\alpha=2 \pi r$.
By (\ref{eqn:Ab2}), it corresponded to the limit $q \to 0$.
In contrast to it, 
here we consider the following scaling limit,
$$
t_{*} \to \infty, \quad \alpha \to \infty \quad
\mbox{with} \quad \frac{t_{*}}{\alpha}=a=\mbox{const.}
$$
By (\ref{eqn:Ab2}), it gives the limit $q \to 1$.
Since
$
\lim_{q \uparrow 1} \sin_q(\pi z)=\sin(\pi z),
$
and $\sin(\pi i z)=i \sinh(\pi z)$, we obtain the limits
\begin{equation}
A_N(x;a) \equiv
\lim_{\substack{t_{*} \to \infty, \alpha \to \infty, \cr t_{*}/\alpha=a}}
A_N^{\alpha}(t_{*}-t,x)
=\frac{1}{2 N a} \coth \left( \frac{x}{2 N a} \right), 
\label{eqn:hyper1}
\end{equation}
\begin{eqnarray}
&& \Phi^A_{\xi, u_k}(z; N, a)
\equiv \lim_{\substack{t_{*} \to \infty, r \to \infty, \cr t_{*}/r=2 \pi a}}
\Phi^A_{\xi, u_k}(z; N, r, t_{*})
\nonumber\\
&& \quad
= \frac{\sinh[(\overline{u}_{\delta}+z-u_k)/2Na]}
{\sinh(\overline{u}_{\delta}/2Na)}
\prod_{\substack{1 \leq \ell \leq N, \cr \ell \not= k}}
\frac{\sinh[(z-u_{\ell})/2Na]}{\sinh [(u_k-u_{\ell})/2Na]},
\,1 \leq k \leq N, \, z \in \C.
\label{eqn:hyper2}
\end{eqnarray}

In this scaling limit the SDE (\ref{eqn:eBES1}) discussed in 
Section \ref{sec:introduction} becomes a temporally homogeneous process on $\R_+$
$$
dX(t)=dB(t)+ \frac{1}{2a} \coth \left(\frac{X(t)}{2a} \right),
\quad t \in [0, \infty).
$$
Similarly (\ref{eqn:hyper1}) and (\ref{eqn:hyper2}) enable us
to discuss the determinantal process
of a hyperbolic version of the Dyson model such that
\begin{equation}
d X^A_j(t) = dB_j(t)
+\frac{1}{2 Na} \sum_{\substack{1 \leq k \leq N, \cr k \not=j}}
\coth \left( \frac{X^A_j(t)-X^A_k(t)}{2Na} \right) dt
+ \frac{1}{2Na} \coth \left( \frac{\overline{X}^A_{\delta}(t)}{2Na} \right) dt,
\label{eqn:hSDE1}
\end{equation}
$1 \leq j \leq N$, $t \geq 0$,
summation of which over $j=1,2, \dots, N$ gives
\begin{equation}
d \overline{X}^A_{\delta}(t)
= \sqrt{N} d B(t)
+\frac{1}{2a} \coth \left( \frac{\overline{X}^A_{\delta}(t)}{2Na} \right) dt.
\quad t \geq 0,
\label{eqn:hSDE2}
\end{equation}
Note that $\lim_{q \uparrow 1} \Gamma_q(z)=\Gamma(z)$
and in this limit (\ref{eqn:gammaq2}) with $z \to iz$ gives 
$\sinh(\pi z)=-\pi i/\{\Gamma(iz) \Gamma(1-iz)\}$, and thus
the above limit is expressed also by using the gamma functions.
If we take the further limit $a \to \infty$ 
with $\delta/a={\rm constant} >0$ in (\ref{eqn:hSDE1}), 
we obtain the Dyson model (\ref{eqn:DysonSDE1}).
In the argument in the present paper,
analyticity of functions $\Phi^A_{\xi, u_k}(z), 1 \leq k \leq N$ plays
an essential role to determine correlation kernels of
determinantal processes.
Since (\ref{eqn:Phib1}) is an entire function of $z$,
the calculations and results for the trigonometric version
of the Dyson model given in Section \ref{sec:reduction_kernel}
will be readily mapped to those for the hyperbolic versions with 
(\ref{eqn:hSDE1}) and (\ref{eqn:hSDE2}) 
by analytic continuation of $z$ and suitable change of 
parameters.

\vskip 0.3cm
\noindent{\bf Acknowledgements} \quad
The present author would like to thank 
Saburo Kakei and Tomoyuki Shirai
for useful comments on the manuscript.
He also thank Masatoshi Fukushima for valuable comments
on Villat's kernel and the Komatu-Loewner evolution.
This work is supported in part by
the Grant-in-Aid for Scientific Research (C)
(No.26400405) of Japan Society for
the Promotion of Science.




\begin{thebibliography}{99}
\bibitem{AAR99}
Andrews, G. E., Askey, R., Roy, R.:
{\it Special Functions},
Cambridge : Cambridge University Press, 1999

\bibitem{BF08}
Bauer, R. O., Friedrich, R. M.:
On chordal and bilateral SLE in multiply connected domain.
Math. Z. {\bf 258}, 241-265 (2008)

\bibitem{BKPV09}
Ben Hough, J., Krishnapur, M., Peres, Y., Vir\'ag, B.:
{\it Zeros of Gaussian Analytic Functions and
Determinantal Point Processes}.
Province R. I. :
Amer. Math. Soc., 2009

\bibitem{Betea11}
Betea, D.:
Elliptically distributed lozenge tilings of a hexagon.
{\sf arXiv:math-ph/1110.4176}

\bibitem{BBO05}
Biane, P., Bougerol, P., O'Connell, N.:
Littelmann paths and Brownian paths.
Duke Math. J. {\bf 130}, 127-167 (2005)

\bibitem{BGR10}
Borodin, A., Gorin, V., Rains, E. M.:
$q$-distributions on boxed plane partitions.
Sel. Math. (N. S.) {\bf 16}, 731-789 (2010)

\bibitem{BR05}
Borodin, A., Rains, E. M.:
Eynard-Mehta theorem, Schur process and their Pfaffian analogs. 
J. Stat. Phys. {\bf 121}, 291-317 (2005) 

\bibitem{Bru91}
Bru, M. F.:
Wishart process.
J. Theoret. Probab. {\bf 3}, 725-751 (1991)

\bibitem{Car03}
Cardy, J.: 
Stochastic Loewner evolution and Dyson's circular ensemble.
J. Phys. A: Math. Gen. {\bf 36}, L379-L386 (2003)

\bibitem{Dys62}
Dyson, F. J. :
A Brownian-motion model for the eigenvalues of a random matrix.
J. Math. Phys. {\bf 3}, 1191-1198 (1962)

\bibitem{Dys72}
Dyson, F. J.:
Missed opportunities.
Bull. Amer. Math. Soc. {\bf 78}, 635-652 (1972)

\bibitem{EM98}
Eynard, B., Mehta, M. L. :
Matrices coupled in a chain: I.
Eigenvalue correlations.
J. Phys. A {\bf 31}, 4449-4456 (1998)

\bibitem{For90a}
Forrester, P. J.:
Exact solution of the lock step model of vicious walkers.
J. Phys. A: Math. Gen. {\bf 23}, 1259-1273 (1990)

\bibitem{For90b}
Forrester, P. J.:
Theta function generalizations of some constant term
identities in the theory of random matrices.
SIAM J. Math. Anal. {\bf 21}, 270-280 (1990)

\bibitem{For10}
Forrester, P. J.:
{\it Log-gases and Random Matrices}.
London Mathematical Society Monographs, Princeton:
Princeton University Press, 2010

\bibitem{Ful04}
Fulmek, M.:
Nonintersecting lattice paths on the cylinder.
S\'eminaire Lotharingien Combin. {\bf 52}, article B52b (2004)

\bibitem{Gos01}
Gosper, R. W.:
Experiments and discoveries in $q$-trigonometry,
In:
{\it Symbolic Computation, Number Theory, Special Functions,
Physics and Combinatorics}.
F. G. Garvan and M. E. H. Ismail (ed),
Dordrecht, Kluwer, 2001, pp.79-105

\bibitem{Gra99}
Grabiner, D. J.:
Brownian motion in a Weyl chamber.
non-colliding particles, and random matrices. 
Ann. Inst. Henri Poincar\'e,
Probab. Stat. {\bf 35}, 177-204 (1999)

\bibitem{Gra02}
Grabiner, D. J.:
Random walk in an alcove of an affine Weyl group,
and non-colliding random walks on an interval.
J. Combin. Theory Ser. A {\bf 97}, 285-306 (2002)

\bibitem{HW96}
Hobson, D. G., Werner, W.:
Non-colliding Brownian motions on the circle.
Bull. London Math. Soc. {\bf 28}, 643-650 (1996)

\bibitem{Joh02}
Johansson, K.:
Non-intersecting paths, random tilings and random matrices. 
Probab. Th. Rel. Fields {\bf 123}, 225-280 (2002)

\bibitem{KM59}
Karlin, S., McGregor, J.:
Coincidence probabilities.
Pacific J. Math. {\bf 9}, 1141-1164 (1959)

\bibitem{Kat13a}
Katori, M.:
Determinantal martingales and noncolliding diffusion processes.
Stochastic Process. Appl. {\bf 124}, 3724-3768 (2014)

\bibitem{Kat13b}
Katori, M.:
Determinantal martingales and correlations of
noncolliding random walks.\\
{\sf arXiv:math.PR/1307.1856}

\bibitem{KT04}
Katori, M., Tanemura, H.:
Symmetry of matrix-valued stochastic processes and
noncolliding diffusion particle systems. 
J. Math. Phys. {\bf 45}, 3058-3085 (2004)

\bibitem{KT10}
Katori, M., Tanemura, H.:
Non-equilibrium dynamics of Dyson's model
with an infinite number of particles.
Commun. Math. Phys. {\bf 293}, 469-497 (2010)

\bibitem{KT11}
Katori, M., Tanemura, H.:
Noncolliding processes, matrix-valued processes
and determinantal processes. 
Sugaku Expositions (AMS) 
{\bf 24}, 263-289 (2011)

\bibitem{KT13}
Katori, M., Tanemura, H.:
Complex Brownian motion representation
of the Dyson model.
Electron. Commun. Probab. {\bf 18}, no.4, 1-16 (2013)

\bibitem{KO01}
K\"onig, W., O'Connell, N.:
Eigenvalues of the Laguerre process as non-colliding 
squared Bessel process.
Elec. Commun. Probab. {\bf 6}, 107-114 (2001)

\bibitem{Kra05}
Krattenthaler, C.:
Advanced determinant calculus: a complement.
Linear Algebra Appl. {bf 411}, 68-166 (2005)

\bibitem{Kra07}
Krattenthaler, C.:
Asymptotics for random walks in alcoves of affine Weyl groups.
S\'eminaire Lotharingien Combin. {\bf 52}, article B52i (2007)

\bibitem{LW13}
Liechty, K., Wang. D.:
Nonintersecting Brownian motions on the unit circle.
{\sf arXiv:math.PR/1312.7390}

\bibitem{Mac72}
Macdonald, I. G. :
Affine root systems and Dedekind's $\eta$-function.
Invent. Math. {\bf 15}, 91-143 (1972)

\bibitem{Meh04}
Mehta, M. L.:
{\it Random Matrices}. 3rd edition, 
Amsterdam: Elsevier, 2004

\bibitem{NF98}
Nagao, T., Forrester, P. J.:
Multilevel dynamical correlation functions for Dyson's
Brownian motion model of random matrices.
Phys. Lett. {\bf A247}, 42-46 (1998)

\bibitem{NF03}
Nagao, T., Forrester,  P.J.:
Dynamical correlations for circular ensembles
of random matrices, 
Nucl. Phys. {\bf B660} 557-578 (2003)

\bibitem{RS06}
Rosengren, H., Schlosser, M.:
Elliptic determinant evaluations and the Macdonald identities
for affine root systems.
Compositio Math. {\bf 142}, 937-961 (2006)

\bibitem{Sch07}
Schlosser, M.:
Elliptic enumeration of nonintersecting lattice paths.
J. Combin. Theory Ser. A {\bf 114}, 505-521 (2007)

\bibitem{ST03}
Shirai, T., Takahashi, Y.: 
Random point fields associated with certain
Fredholm determinants I:
fermion, Poisson and boson point process.
J. Funct. Anal.
{\bf 205}, 414-463 (2003)

\bibitem{Sos00}
Soshnikov, A. : 
Determinantal random point fields.
Russian Math. Surveys {\bf 55}, 923-975 (2000)

\bibitem{Spo87}
Spohn, H.:
Interacting Brownian particles:
a study of Dyson's model.
In:
{\it Hydrodynamic Behavior and Interacting Particle Systems}. 
G. Papanicolaou (ed),  
IMA Volumes in Mathematics and its Applications, {\bf 9}, Berlin: 
Springer-Verlag, 1987, pp.151-179

\bibitem{TV97}
Tarasov, V., Varchenko, A.: 
Geometry of $q$-hypergeometric functions,
quantum affine algebras and elliptic quantum groups.
Ast\'erisque {\bf 246}, 1-135 (1997)

\bibitem{War02}
Warnaar, S. O.:
Summation and transformation formulas
for elliptic hypergeometric series.
Constr. Approx. {\bf 18}, 479-502 (2002)

\bibitem{WW27} 
Whittaker, E. T., Watson, G. N.:
{\it A Course of Modern Analysis}. 4th edition,
Cambridge: Cambridge University Press, 1927

\bibitem{Zhan04}
Zhan, D.:
Stochastic Loewner evolution in doubly connected domains.
Probab. Theory Relat. Fields {\bf 129}, 340-380 (2004)
\end{thebibliography}
\end{document}